\documentclass{article}

\title{Parameter Estimation for the Stochastically Perturbed Navier-Stokes Equations\footnote{to appear in Stochastic Processes and their Applications. 
\href{http://dx.doi.org/10.1016/j.spa.2010.12.007}{DOI:10.1016/j.spa.2010.12.007} }}
\author{Igor Cialenco
\footnote{Department of Applied Mathematics,
Illinois    Institute of Technology, Email: \url{igor@math.iit.edu}} \, and
  Nathan Glatt-Holtz\footnote{Department of Mathematics, Indiana University, Bloomington. Email: \url{negh@indiana.edu}}\\
 }
\usepackage{amsfonts, amssymb, amsmath, amsthm}
\usepackage[colorlinks=true, pdfstartview=FitV, linkcolor=blue,
            citecolor=blue, urlcolor=blue]{hyperref}
\usepackage[usenames]{color}
\usepackage{accents}
\usepackage{comment}

\usepackage{url}
\makeatletter\def\url@leostyle{%
  \@ifundefined{selectfont}{\def\UrlFont{\sf}}{\def\UrlFont{\scriptsize\ttfamily}}} \makeatother

\numberwithin{equation}{section}

\newtheorem{Thm}{Theorem}[section]
\newtheorem{Lem}[Thm]{Lemma}
\newtheorem{Prop}[Thm]{Proposition}
\theoremstyle{definition}
\newtheorem{Cor}[Thm]{Corollary}

\newtheorem{Rmk}[Thm]{Remark}

\newtheorem*{Rmk*}{Remark}

\newcommand{\pd}[1]{\partial_{#1}}

\newcommand{\E}{\mathbb{E}}
\newcommand{\Prb}{\mathbb{P}}

\includecomment{commentCoAuthor}
\includecomment{commentToDo}

\begin{document}
\markboth{I. Cialenco, N. Glatt-Holtz}
{Parameter Estimation for the Stochastic Navier-Stokes Equations}

\maketitle

\begin{abstract}
 We consider a parameter estimation problem to determine  the viscosity $\nu$ of a stochastically perturbed 2D Navier-Stokes system. We derive several different classes of   estimators based on the first $N$ Fourier modes of a  single sample path observed on a finite time interval. We study the consistency and asymptotic normality of these estimators.  Our analysis treats strong, pathwise solutions for both the periodic and bounded domain cases in the presence of an additive white (in time) noise.\\
\end{abstract}

{\noindent \small
}

\section{Introduction}

The theory of stochastic partial differential equations (SPDEs) is a
rapidly developing field of pure and applied mathematics.  These equations
are used to describe the evolution of dynamical systems in the presence of persistent
spatial-temporal uncertainties.   When considering nonlinear processes one encounters
many new, fundamental and mathematically challenging problems for SPDEs, with
important applications in physics and applied sciences.

While the general form of a particular SPDE
is commonly derived from the fundamental properties
of the underlying processes under study, frequently parameters
arise in the formulation which need to
be specified or determined on the basis of some sort of empirical observation.
In such situations, the so called problem of \emph{parameter estimation}
arises naturally: \emph{under the assumption that a
phenomenon of interest follows the dynamics of an SPDE, and given that some realizations
of this process are measured, we wish to find the unknown parameters appearing in the model,
such that the equations fit or predict as much as possible the observed data. }

Actually, the development of methods to estimate parameters appearing in a model
serve practical considerations for two reasons.
On the one hand we may be confident in the model, but have an incomplete
knowledge of the physical parameters appearing therein.  An ``estimator'' of the true parameter
therefore provides a means to measure these unknowns.  On the other
hand, we may already possess accurate knowledge of the physical
quantities involved in the model, but lack confidence in the validity of the
underlying model.  In this situation finding an ``estimator'' will be
the first step in testing and validating the model.

Since the solution of an SPDE is a random variable, this inverse problem of
finding the true parameters is treated by methods from stochastic analysis and statistics.
In this work we will follow a continuous time
approach and assume that the solution $U = U_\nu(t,\omega)$ of the SPDE
is observed for every time $t$ over an interval $[0,T]$.
We note that different types of methods and approaches are used to study
inverse problems for deterministic PDEs, and we refer the reader to
\cite{Isakov1, Kirsch1} and references therein.

A core notion in the theory of statistical inference for stochastic processes
is the so called `regularity' of the family of probability measures associated
to the set of possible values $\Theta$ of the parameter of interest $\nu$.
Note that $\nu$ could be a vector in general.   Let
$H$ be the function space where the solution evolves and for each $\nu \in \Theta$
denote by $\mathbb{P}_\nu^T$, the probability
measures on $C([0,T]; H)$ generated by the solution $U_\nu(t), \ 0\leq t\leq T$.
We say that a model is `regular' if any two probability measures
from the family $\{ \mathbb{P}_\nu^T, \ \nu\in\Theta\}$, are mutually absolutely continuous.
On the other hand the model is said to be `singular' if these measures are mutually singular.

For regular models one approach to the parameter estimation problem is
to consider the Maximum Likelihood Estimator (MLE) $\widehat{\nu}$ of $\nu$.
This type of estimator is obtained by fixing a reference value $\nu_{0}$ and
then maximizing the Radon-Nikodym derivative
or Likelihood Ratio $d\mathbb{P}_\nu^T/d\mathbb{P}_{\nu_0}^T$ with respect to $\nu$.
Usually $\widehat{\nu}\neq \nu$ and the problem is to study the convergence
of these estimators to the true parameter as more information arrives (for example
as time passes or by decreasing the amplitude of the noise).
In contrast, each singular model requires an individual approach, and usually the true
parameter can be found exactly, without any limiting procedure (at least if the solution
is observed continuously.)

Statistical inference for finite dimensional systems of stochastic
differential equations (SDEs) have been studied widely and
provide instructive examples of both `regular' and `singular' problems.
Typically estimating the drift coefficient for an SDE is a
regular problem which may be treated with an MLE.
Here the likelihood ratio can be determined by Girsanov type theorems.
By contrast, estimating the diffusion coefficient is a singular problem and
in this case one can find the diffusion coefficient by measuring the
quadratic variation of the process.
In general there exist necessary and sufficient conditions for the regularity
for (finite dimensional) SDEs.  See the monographs \cite{KutoyantsBook2004},
\cite{LiptserShiryayev}, and references therein for a comprehensive treatment.

It turns out that the parameter estimation problem for infinite
dimensional systems (SPDEs) is, in many cases,  a singular problem where
one can find the parameter  $\nu$ ``exactly'' on any finite interval of time.
In particular this has been shown in the case of linear stochastic
parabolic equations with the parameter of interest in the drift appearing next to
the highest order differential operator.
Note that this is in direct contrast to most of the
corresponding finite dimensional processes where one has to
observe a sample path over an infinite time horizon or to decrease
the amplitude of the noise term in order to get similar results.
One of the first significant works in the theory of statistical inference
for SPDEs that explorers this singularity is \cite{HubnerRozovskiiKhasminskii}.
The idea in this work is to approximate the original singular problem by a
sequence of regular problems for which MLEs exist.
This approximation is carried out by considering Galerkin-type projections
of the solution onto a finite-dimensional space where the estimation
problem becomes regular.
They prove that as the dimension of the projection increases
the corresponding MLEs converge to the true parameter.
In \cite{HuebnerLototskyRozovskii97,
HuebnerRozovskii, LototskyRozovskii1999, LototskyRozovskii2000},
the problem has been extended to a
general class of linear parabolic SPDE driven by additive noise
and the convergence of the estimators
has been classified in terms of the
order of the corresponding differential operators.
%For equations driven by multiplicative noise see \cite{IgorMultiplicativeFBM2010, IgorSergey2007}.
For recent developments and  other types of inference
problems for \emph{linear SPDEs} see the survey paper
 \cite{Lototsky2009Survey} and containing references.

While the linear theory has been extensively studied
in the framework described above it seems that,
to the best of our knowledge, no similar results
have been established for \emph{nonlinear SPDEs}.
We therefore embark in this and concurrent work
\cite{IgorNathanMultiplicativeNS2010}
on a study of parameter estimation problems for certain fundamental
nonlinear SPDEs from fluid dynamics.

Note that for the linear case, key properties
such as efficiency and asymptotic normality of the estimators,
are proven by making essential use of the exact long time behavior of the moments
of the Fourier coefficients of the solutions.
In the case of nonlinear equations, for example stochastic
equations from mathematical fluid dynamics, the problem is much more delicate, due to
the (highly nontrivial) coupling of the Fourier modes.

From the point of view of applications this work is
motivated in particular by recent developments
in the area of Geophysical Fluid Dynamics (GFD)
where the theory of SPDEs is now playing an important role.
See, for example, \cite{Penland1, PenlandEwald1, PenlandMatrosova1,
EwaldPetcuTemam,  GlattHoltzZiane1,
GlattHoltzTemam2, DebusscheGlattHoltzTemam1}.
For this developing field, novel `inverse' methods are clearly needed.
While the problems we consider initially are toy
models in comparison to large scale circulation models
such as the Primitive Equations,
we are optimistic that the methods and insights developed for
simple nonlinear SPDEs will eventually serve the wider goal of
extending our understanding to a more physically realistic setting.

In this work we consider the $2D$
Navier-Stokes equations forced with an additive white noise:
\begin{subequations}\label{eq:NSESTochastic}
 \begin{align}
   dU + ((U \cdot \nabla) U - \nu \Delta U  + \nabla P) dt &=
       \sigma dW,
     \label{eq:NSEMomentum}\\
  \nabla \cdot U &= 0,
     \label{eq:NSEDivFreeCond}\\
    U(0) &= U_0,
     \label{eq:NSEInitData}
 \end{align}
\end{subequations}
which describe the flow of a viscous,
incompressible fluid. Here $U = (U_1, U_2)$ and  $P$ respectively
represent the
velocity field and the pressure.   The coefficient $\nu > 0$ corresponds to
the kinematic viscosity of the fluid, and it will be the parameter of interest.
The goal of our analysis will be to find a suitable
estimator $\hat{\nu} = \hat{\nu}(\omega)$ which is a functional of
a single sample path $U(\omega)$ observed over a finite and
fixed time interval $[0,T]$.

We assume that the governing equations (\ref{eq:NSESTochastic})
evolve over a domain  $\mathcal{D}$.   Throughout this work
we will consider two possible boundary conditions.  On the one hand
we may suppose that the flow occurs over all of $\mathbb{R}^2$, take
$\mathcal{D} = [-L/2, L/2]^{2}$ for some $L > 0$ and prescribe the
\emph{periodic boundary condition}:
\begin{equation}\label{eq:PeriodicBC}
 U(\mathbf{x} + L \mathbf{e}_j, t) = U(\mathbf{x},t), \
 \textrm{ for all } \mathbf{x} \in \mathbb{R}^2, t \geq 0;
 \quad \int_{\mathcal{D}} U(\mathbf{x}) d\mathbf{x} = \mathbf{0} \, .
 \footnote{This second condition may be added with no loss of
 generality and slightly simplifies the analysis.  See for instance \cite{Temam4}.}
\end{equation}
We also consider the case when $\mathcal{D}$ is a bounded subset of
$\mathbb{R}^2$ with a smooth boundary $\partial\mathcal{D}$ and assume
the \emph{Dirichlet (no slip) boundary condition}:
\begin{equation}\label{eq:BCDirchlet}
  U(\mathbf{x}, t) = 0 \textrm{ for all }
  \mathbf{x} \in \partial\mathcal{D}, t \geq 0.
\end{equation}

\emph{The stochastic forcing} we consider is an additive
space-time noise colored in space. Formally, we may write
\begin{equation}  \label{eq:noiseStructure}
  \sigma dW = \sum_k \lambda_k^{-\gamma} \Phi_k dW_{k},
\end{equation}
where $\Phi_k$ are the eigenfunctions of the Stokes operator,
$\lambda_k$ represent the associated
eigenvalues, and $W_k,\ k\geq 1,$ are one dimensional independent Brownian motions.
We assume that $\gamma$ is a real parameter greater than $1$ which
guarantees some spatial smoothness in the forcing.
We may also formally derive (see e.g. \cite{ZabczykDaPrato1}) the space-time
correlation structure of the noise term
\begin{equation*}
  \E( \sigma dW (\mathbf{x}, t) \sigma dW (\mathbf{y}, s))
  	= K(\mathbf{x},\mathbf{y}) \delta_{t -s}\, ,
\end{equation*}
where $K(\mathbf{x},\mathbf{y}) = \sum_{k \geq 1} \lambda_{k}^{-2\gamma}
\Phi_{k}(\mathbf{x})\Phi_{k}(\mathbf{y})$.

We should mention that the Stochastic Navier-Stokes equations
in both $2$ and $3$ dimensions and under much more general stochastic forcing
conditions have been extensively studied.   See, for instance,
\cite{BensoussanTemam, Breckner, BrzezniakPeszat, Cruzeiro1,
FlandoliGatarek1, GlattHoltzZiane2, MikuleviciusRozovskii2} and
containing references.

Since the parameter of interest $\nu$ appears next
to the highest order differential operator,
the linear analogue of \eqref{eq:NSESTochastic} is
singular as we described above.  With this in mind
we expected that the full nonlinear model might also be singular.
As developed below, $\nu$ may be found exactly from
a single observation over a finite time window which
suggests that this singular structure is preserved in
this nonlinear case.

The starting point of our analysis, the derivation of an estimator for $\nu$,
follows methods already
developed for the linear case (see references mentioned above).
We project \eqref{eq:NSESTochastic} down to a finite
dimensional space, and for each $N$ we arrive at a system of the form
\begin{displaymath}
   dU^N + (\nu AU^N + P_NB(U))dt = P_N\sigma dW, \quad U(0) = U_0,
\end{displaymath}
where $P_N$ is the projection operator on the finite dimensional space generated by the
first $N$ Fourier eigenvalues of the Stokes operator.
We then formally compute the MLEs associated to these systems,
and take them as an ansatz for our estimators.
In the course of the analysis we  introduce an
additional degree of freedom, a parameter $\alpha$, which we may
carefully tune to compensate for the nonlinear term.
We arrive finally at the following three classes of estimators:
\begin{equation}\label{eq:estInitialForm}
\begin{split}
\widetilde{\nu}_N & =
 -\frac{\int_0^T \langle A^{1 + 2 \alpha} U^N, dU^N \rangle
   + \int_0^T \langle A^{1 + 2 \alpha}U^N, P_N B(U) \rangle dt  }
 { \int_0^T |A^{1+ \alpha}U^N|^2 dt}, \\
\check{\nu}_N & = -  \frac{\int_0^T \langle A^{1 + 2 \alpha} U^N, dU^N \rangle
   + \int_0^T \langle A^{1 + 2 \alpha}U^N, P_N B(U^N) \rangle dt  }
 { \int_0^T |A^{1+ \alpha}U^N|^2 dt}, \\
 \hat{\nu}_N & =  - \frac{\sum\limits_{k =1}^{N} \lambda_k^{1 + 2 \alpha}
           (u_k^2(T) - u_k^2(0) - T \lambda_k^{-2\gamma}) }
        {2 \sum\limits_{k =1}^{N}\lambda_k^{2 + 2 \alpha} \int_0^T u_k^2 dt}.
\end{split}
\end{equation}
Here $u_k := (U, \Phi_k)$ represents the $k$th (generalized) Fourier
mode of the solution $U$.

The main result in this work establishes the following properties for
the proposed estimators:

\begin{Thm}\label{thm:MainResultBasic}
  Suppose that $U = U(\omega)$ is a single sample path solution of
  (\ref{eq:NSESTochastic}), (\ref{eq:PeriodicBC}) or
  (\ref{eq:NSESTochastic}), (\ref{eq:BCDirchlet}) observed on a finite interval of time $[0,T]$.
  Assume that   (\ref{eq:NSESTochastic}) is forced with a white noise process of the
  form (\ref{eq:noiseStructure}) where $\gamma > 1$,\footnote{In the case  (\ref{eq:BCDirchlet})
            we assume, for technical reasons, an upper bound on $\gamma$,
            $\gamma < 1 + 1/4$ as well.  See below.}
 and suppose that $\alpha > \gamma - 1$.
Then, given a suitably regular initial vector
field $U_0$,
\begin{itemize}
\item[(i)] the functionals $\widetilde{\nu}_N,  \check{\nu}_N, \widehat{\nu}_N$
 defined by \eqref{eq:estInitialForm}
 are weakly consistent estimators of the parameter $\nu$, i.e.
\begin{equation*}%\label{eq:consistencynuHat}
     \lim_{N \rightarrow \infty} \widetilde{\nu}_N  = \lim_{N \rightarrow \infty} \check{\nu}_N = \lim_{N \rightarrow \infty} \widehat{\nu}_N = \nu
\end{equation*}
in probability.
\item[(ii)] if we assume further that $\alpha > \gamma - 1/2$, then
    $\widetilde{\nu}_{N}$ is
    asymptotically normal with rate $N$ i.e.
    \begin{displaymath}%{equation}\label{eq:asymtoticNormalnuHat}
     N (\widetilde{\nu}_N  -\nu)
        \overset{d}{\longrightarrow} \eta
  \end{displaymath}%{equation}
    (converges in distribution) where $\eta$ is a mean zero, normally
    distributed random variable.
  \end{itemize}
\end{Thm}

While we are able to prove the strongest convergence
results for $\widetilde{\nu}_{N}$, this estimator is intractable numerically
and even analytically.  This is because $\widetilde{\nu}_{N}$
depends on all of the
Fourier modes of the solution in a highly nonlinear fashion.
At the other extreme is $\hat{\nu}_{N}$ which is much more
straightforward to compute but is expected to have a slower
rate of convergence to the actual parameter $\nu$.  The estimator
$\check{\nu}_N$ is a compromise between the two extremes
since it depends only on the knowledge of the first $N$ eigenmodes
but retains some of the complex structure of the nonlinear term.
Although at the present time we are not able to prove this,
we expect that $\check{\nu}_N$ has a faster rate of convergence
than $\hat{\nu}_{N}$.
We  conjecture, in Section \ref{sec:asymtotic-normality},
that $\check{\nu}_N$ is also asymptotically normal with the same variance and
rate of convergence as $\widetilde{\nu}_N$.
Given the explicit formulas for the estimators, \eqref{eq:estInitialForm},
all these questions, including the effect of the free parameter $\alpha$
on the rate of convergence, can be studied by means of numerical simulations,
which the authors plan to undertake in a separate forthcoming paper.

While the form of the proposed estimators
and the general statements of the main results in this work are similar
to previous works in the linear case, fundamental new difficulties
arise  which require one to take a novel approach
for the analysis. This is of course due to the complex structure
of the nonlinear term appearing in (\ref{eq:NSESTochastic}) which
couples, in an intricate way, all of the modes $u_{k} = (U, \Phi_{k})$.
In contrast to the linear case, we lose for example any explicit spectral
information about the elements $u_{k}$.  This coupling also means
that the $u_{k}$ are not expected to be independent.

To overcome these difficulties the analysis relies on a careful decomposition
of the solution $U = \bar{U} + R$.  Here $\bar{U}$ satisfies a linear system
where the modes are independent.  Crucially, a complete spectral picture is
obtainable for $\bar{U}$.  On the other hand, $R$, while depending in a
complicated way on the full solution $U$ is more regular in comparison
to $\bar{U}$.  This is because $R$ is not directly forced by the noise terms
$\sigma dW$.  For this point the analysis, particularly in the case of bounded
domains, requires a delicate treatment of the nonlinear term.

Due to these technical issues, we were able to establish asymptotic normality
only for $\widetilde{\nu}_N$.
It is interesting that $\widehat{\nu}_N$ is a
consistent estimator for $\nu$ and it is
the same as the MLE of the corresponding linear equation
(the stochastic Stokes equation).
This effect can be explained as follows: since the nonlinear term
$B(U)$ is in some sense `lower order'  it fails to destroy the information
about $\nu$; $\nu$ remains observable in a similar manner to the linear case.

The exposition of the paper is organized as follows: In
Section~\ref{sec:MathSetting} we lay the theoretical foundations for
this work reviewing the relevant mathematical theory for the
stochastic Navier-Stokes equations. We establish some
crucial spectral information concerning the linear system associated
to (\ref{eq:NSESTochastic}). We also recall in this section some
particular variants on the law of large numbers and the central limit
theorem. Section~\ref{sec:DerivationOfTheEstimator} sketches
the derivation of the estimators $\widetilde{\nu}, \ \check{\nu}, \ \widehat{\nu}_N$.
We conclude the section with a strict formulation of the main results.
The proof of the main theorem is carried out in
Section~\ref{sec:ProofOfMainRes} in a series of modular substeps.
We first study the regularity of the `residual' $R$ that appears after we
`subtract off' the noise term appearing in (\ref{eq:NSESTochastic})
via the linear Stokes equation.
As an immediate application we are able to determine some precise rates
for the denominators appearing in the estimators \eqref{eq:estInitialForm} . Using these
rates we successively analyze the consistency of the estimators.
The final subsection treats the question of asymptotic normality with
the help of a central limit theorem for martingales.

\section{Mathematical Setting of the Problem}
\label{sec:MathSetting}

We begin by recalling the mathematical background for the stochastic
Navier-Stokes Equations and then review some general results from
probability theory that will be used in the sequel.

\subsection{The Stochastic Navier-Stokes Equation}
We first describe how (\ref{eq:NSESTochastic}) is recast as an
infinite dimensional stochastic evolution equation of the form
\begin{equation}\label{eq:SNSEabs}
  \begin{split}
    dU + (\nu AU + B(U))dt
     &= \sigma dW\, , \\
    U(0) &= U_0 .
 \end{split}
\end{equation}
The basic functions spaces are designed to capture both the
boundary conditions and the divergence free nature of the flow.

We first consider the spaces associated with a \emph{Dirichlet boundary
condition} (\ref{eq:BCDirchlet}).  Let
$H := \{ U \in  L^2(\mathcal{D})^2: \nabla \cdot U = 0,  U \cdot n = 0\}$,
where $n$ is the outer pointing unit normal to $\partial \mathcal{D}$.
$H$ is endowed as a Hilbert space with the $L^2$ inner product
$(U^{\flat}, U^\sharp) = \int_{\mathcal{D}}  U^{\flat} U^\sharp dx$ and
associated norm $| U | = ( U, U)^{1/2}$.
The Leray-Hopf projector, $P_H$, is defined as the orthogonal projection
of $L^2(\mathcal{D})^d$ onto $H$.  We next take
$V := \{ U \in H^1_0(\mathcal{D})^2 : \nabla \cdot U = 0 \}$
and endow this space with the inner product
$((U^{\flat}, U^\sharp)) = \int_{\mathcal{M}} \nabla
U^{\flat} \cdot \nabla U^\sharp
  d\mathcal{M}$.
  Due to the Dirichlet boundary condition, (\ref{eq:BCDirchlet}),
the Poincar\'{e} inequality $|U|
\leq c\|U\|$ holds for $U \in V$ justifying
this definition.

The definitions for $H$ and $V$ are slightly different for the case
of \emph{periodic boundary conditions} (\ref{eq:PeriodicBC}). We
take $\mathcal{D} = [-L/2, L/2]^2$ and define the spaces
$L^2_{per}(\mathcal{D})^2$, $H^1_{per}(\mathcal{D})^2$ to be
the families of vector fields $U = U(\mathbf{x})$ which are $L$
periodic in each direction and which belong respectively to
$L^2(\mathcal{O})^2$ and $H^1(\mathcal{O})^2$ for every open bounded
set $\mathcal{O} \subset \mathbb{R}^2$.  We now define
\begin{displaymath}
  H = \left\{ U \in L^2_{per}(\mathcal{D})^2:
    \nabla \cdot U = 0,
    \int_{\mathcal{D}} U(\mathbf{x}) d\mathbf{x} = \mathbf{0}
  \right\},
\end{displaymath}
and
\begin{displaymath}
  V = \left\{ U \in H^1_{per}(\mathcal{D})^2:
    \nabla \cdot U = 0,
    \int_{\mathcal{D}} U(\mathbf{x}) d\mathbf{x} = \mathbf{0}
  \right\}.
\end{displaymath}
$H$ and $V$ are endowed with the norms $| \cdot |$
and $\| \cdot \|$ as above.
Note that we impose the mean zero condition for $H$ and $V$
so that the Poincar\'{e} inequality holds.  As
mentioned in the introduction, there is no loss of generality
in imposing this extra assumption.  See, e.g. \cite{Temam4}.

The linear portion of \eqref{eq:NSESTochastic} is captured in the Stokes
operator $A = - P_H \Delta$, which is an unbounded operator from $H$ to $H$
 with the domain $D(A) = H^2(\mathcal{M}) \cap V$.  Since $A$ is self adjoint, with
a compact inverse $A^{-1} : H \rightarrow D(A)$, we may apply the standard
theory of compact, symmetric operators to guarantee the existence of an
orthonormal basis $\{\Phi_k\}_{k \geq 1}$ for $H$ of eigenfunctions of $A$ with the
associated eigenvalues $\{\lambda_k\}_{k \geq 0}$ forming an unbounded, increasing, sequence.
Moreover,
\begin{equation}\label{eq:LambdaOrder}
  \lambda_k \approx \lambda_1 k,
\end{equation}
where the notation $a_n \approx b_n$ means that $\lim_{n\to\infty}a_n/b_n=1$. Also, we will write
$a_n\sim b_n$ when there exists a finite, nonzero constant $c$ such that
$\lim_{n\to\infty}a_n/b_n=c$. For more details about asymptotical behavior of $\{\lambda_k\}_{k\geq 1}$
see for instance \cite{Babenko1, Metivier1} for the no-slip case
(\ref{eq:BCDirchlet}), and \cite{ConstantinFoias1} for the spatially
periodic case (\ref{eq:PeriodicBC}).   Define
$H_N = \textrm{Span} \{\Phi_1, \ldots, \Phi_N\}$,
and take $P_N$ to be the projection from $H$ onto this space.  We
let $Q_N := I - P_N$.

The analysis below relies extensively on the fractional powers of $A$.
Given $\alpha > 0,$ take
$D(A^\alpha) = \left\{
   U \in H:  \sum_k \lambda_k^{2\alpha} |u_k|^2 < \infty \right\}$,
where $u_k = (U, \Phi_k)$.  On this set we may define $A^\alpha$
according to
$A^\alpha U  = \sum_{k} \lambda_k^{\alpha} u_k \Phi_{k}$,
 for $U = \sum_{k} u_{k} \Phi_{k}$.
Classically we have the generalized Poincar\'{e} and inverse Poincar\'{e} estimates
\begin{equation}\label{eq:decompEstimates}
 \begin{split}
    |A^{\alpha_2} P_N U | \leq \lambda^{\alpha_2 - \alpha_1}_N |A^{\alpha_1} P_N U|,
    \quad
    |A^{\alpha_1} Q_N U | \leq
     \frac{1}{\lambda^{\alpha_2 - \alpha_1}_N} |A^{\alpha_2} Q_N U|,\\
 \end{split}
\end{equation}
for any $\alpha_1 < \alpha_2$.

We next describe the stochastic terms in \eqref{eq:NSESTochastic}.
Fix a stochastic basis $\mathcal{S} := (\Omega, \mathcal{F},
\{\mathcal{F}_t\}_{t \geq 0}, \mathbb{P}, \{W_k\}_{k \geq 1})$, that
is a filtered probability space with $\{W_k\}_{k \geq 1}$ a sequence
of independent standard Brownian motions relative to filtration $\{\mathcal{F}_t\}_{t\geq 0}$.
In order to avoid unnecessary complications below we may assume that
$\mathcal{F}_t$ is complete and right continuous (see \cite{ZabczykDaPrato1} for more details).
Writing formally $W = \sum_{k \geq 0} \Phi_k W_k$, $W$ may be viewed as
a cylindrical Brownian motion on $H$.

We briefly recall the classical formalism for infinite-dimensional Wiener process
as in \cite{ZabczykDaPrato1},
\cite{PrevotRockner}.
Consider the collection of Hilbert-Schmidt operators mapping $H$
into $D(A^\beta)$, $\beta \geq 0$.  We denote this family by $L_2(H,D(A^{\beta}))$.  Throughout
this work we assume that $\sigma$, understood as an operator, has the form
\begin{equation}\label{eq:sigmaOperation}
  \sigma \Phi_k = \lambda_k^{-\gamma} \Phi_k.
\end{equation}
We will write
$\sigma dW(t) = \sum_{k\geq 1} \lambda_k^{-\gamma} \Phi_k dW_{k}(t), \ t\geq 0.$
One may check that, for every $\epsilon > 0$,
$\sigma \in L_2(H, D(A^{\gamma - 1/2 - \epsilon}))$.
In particular, given the standing assumption that
$\gamma > 1$, we have $\sigma \in L_2(H,D(A^{1/2}))$.

\subsection{The Stochastic Stokes Equation and Limit Theorems}\label{sec:StokesEqAndLimitThrms}
We next consider the linear system
associated to (\ref{eq:SNSEabs}), which we write in the abstract form:
\begin{equation}\label{eq:StokesStochasticAbs}
   d\bar{U} + \nu A\bar{U} dt
      	= \sum_k \lambda_k^{-\gamma} \Phi_k dW_k, \quad \bar{U}(0) = \bar{U}_0.
\end{equation}
For the purposes  here this system can be analyzed as 2D stochastic heat equation driven by an additive cylindrical Brownian motion (for general results we refer readers to \cite{ZabczykDaPrato1, RozovskiiBook}.)
%Such infinite dimensional linear systems are well studied, see e.g. \cite{ZabczykDaPrato1, RozovskiiBook}.

Let us denote by $\bar{u}_k,\ k\geq 1$, the Fourier coefficients  of the solution $\bar{U}$ with respect to the system $\{ \Phi_k\}_{k}$ in $H$,
i.e. $\bar{u}_k = (\bar{U}, \Phi_k), \ k\geq 1$.
By (\ref{eq:StokesStochasticAbs}), we note that each Fourier mode $\bar{u}_k$
represents a one dimensional stable Ornstein-Uhlenbeck process with dynamics
\begin{equation}\label{eq:stokesSpecEqns}
   d\bar{u}_k + \nu \lambda_k \bar{u}_k dt =  \lambda_k^{-\gamma} dW_k,
   \quad \bar{u}_k(0) = \bar{u}_{0k}, \ k\geq 1.
\end{equation}
It follows from \eqref{eq:stokesSpecEqns} that
\begin{equation}\label{eq:OUSolution}
\bar{u}_k(t) = \bar{u}_k(0)e^{-\nu\lambda_k t} + \lambda_k^{-\gamma}\int_0^t e^{-\nu\lambda_k(t-s)} dW_k(t), \quad k\geq 1, \ t\geq 0.
\end{equation}

In what follows we will use the following auxiliary results about asymptotics  of
the first moments of the Fourier modes  $\bar{u}_k, \ k\geq 1$ (see also Theorem 2.1
in \cite{Lototsky2009Survey}.)
\begin{Lem}\label{thm:RegStokesEqn}
  Suppose that $\bar{U}$ is a solution of (\ref{eq:StokesStochasticAbs}) and
  let $\bar{U}^N := P_N \bar{U}, \ N\geq 1$.
  \begin{itemize}
  \item[(i)]  Assume that $\gamma' < \gamma$ and that $\E |A^{\gamma' -1/2} \bar{U}_0|^{2} < \infty$. Then
    \begin{equation}\label{eq:LinearEquationReg}
      \bar{U} \in
      L^{2}(\Omega; L^{2}_{loc}([0,\infty); D(A^{\gamma'})))\cap
      L^{2}(\Omega; C([0,\infty); D(A^{\gamma'-1/2}))) \, .
    \end{equation}
  \item[(ii)] Suppose that $\bar{U}_0 = 0$, then:
    \begin{equation}    \label{eq:specL2ubark}
      \E \int_0^T \bar{u}_k^2 dt \approx \frac{T \lambda_k^{-(1 + 2 \gamma)}}{2\nu}
      \approx \frac{T \lambda_1^{-(1+2\gamma)} }{2 \nu} k^{-(1 + 2\gamma)}\, ,
    \end{equation}
    and
    \begin{equation} \label{eq:specL4Typeubark}
      \mathrm{Var} \left[\int_0^T \bar{u}_k^2 dt \right]
      \sim \lambda_k^{-(3 + 4 \gamma)}
      \sim k^{-(3 + 4 \gamma)}.
    \end{equation}
\item[(iii)]     Moreover, for $\beta > \gamma$,
   \begin{equation}\label{eq:BlowUpRateLinearSystem}
     \E \int_0^T |A^\beta \bar{U}^N|^2 dt
     \approx \frac{T \lambda_1^{2\beta-2\gamma-1}}{2\nu(2\beta-2\gamma)} \, N^{2\beta-2\gamma}.
      % \lambda_N^{2\beta - 2 \gamma} \sim  N^{2\beta - 2 \gamma}.
   \end{equation}
\end{itemize}
\end{Lem}
\begin{proof}
  The first item is classical and may, for example, be justified with
a Galerkin Scheme or other suitable techniques from the general theory
of existence and uniqueness of the solutions
for stochastic parabolic equations.  See e.g. \cite{ZabczykDaPrato1, RozovskiiBook}.
Using \eqref{eq:OUSolution}, (ii) follows by direct computations of the corresponding moments, and for
the final item we deduce
\begin{align*}\label{eq:RegEstBarU}
      \E \int_{0}^{T} |A^{\beta} \bar{U}^N|^{2}dt
      & =   \E \int_{0}^{T} |\sum_{k=1}^N \lambda_{k}^{\beta} \bar{u}_{k} \Phi_{k}|^{2}dt
      =   \sum_{k=1}^N \lambda_{k}^{2\beta}  \E \int_{0}^{T} \bar{u}_{k}^{2} dt  \\
          &  \approx    \frac{T}{2\nu} \sum_{k=1}^N \lambda_{k}^{2\beta - 1 - 2 \gamma}
            \approx  \frac{T\lambda_1^{2\beta-2\gamma-1}}{2\gamma} \frac{N^{2\beta-2\gamma}}{2\beta-2\gamma},
             % \lambda_N^{2\beta - 2\gamma} \sim N^{2\beta - 2\gamma},
  \end{align*}
where we have made use of (ii), (\ref{eq:LambdaOrder}) in conjunction
  with
  \begin{equation}\label{eq:partialSum}
    \sum_{k=1}^N k^{a} \approx \frac{N^{1+ a}}{a+1}, \quad a > -1.
  \end{equation}
  The proof is complete.
\end{proof}

We finally recall some particular versions of the Law of Large Numbers (LLN)
and the Central Limit Theorem (CLT) which are used to prove consistency and
asymptotic normality of the class of estimators given by \eqref{eq:estInitialForm}.

\begin{Lem}[The Law of Large Numbers]\label{thm:LLN}
Let $\xi_n, \ n\geq 1$, be a sequence of random variables and $b_n, \ n\geq 1$,
an increasing sequence of positive numbers such that $\lim_{n\to\infty}b_n=+\infty$, and
\begin{equation}\label{eq:bndCondLLN}
   \sum\limits_{n=1}^\infty \frac{ \mathrm{Var}\xi_n}{b_n^2} <\infty.
\end{equation}
\begin{itemize}
\item[(i)] If we assume that the random variables $\xi_n, \ n\geq 1$, are independent then
  \begin{equation*}%\label{eq:stongConverLNN}
   \lim\limits_{n\to\infty} \frac{ \sum\limits_{k=1}^n (\xi_k - \mathbb{E}\xi_k)}{b_n} = 0
   \quad a.s.
 \end{equation*}
\item[(ii)] If we suppose only that $\xi_n, \ n\geq 1$, are merely
  uncorrelated random variables, then
  \begin{equation}\label{eq:stongConverLNN}
   \lim\limits_{n\to\infty} \frac{ \sum\limits_{k=1}^n (\xi_k - \mathbb{E}\xi_k)}{b_n} = 0,
 \end{equation}
in probability.
\end{itemize}
\end{Lem}

  \begin{proof}
  See, for example, Shiryaev \cite[Theorem IV.3.2]{ShiryaevBookProbability} for the proof of (i).
Part two, similar to the proof of  Weak LLN, follows from the Markov inequality.
For a fixed $\epsilon > 0$, and for all pairs $m < n$,
we have
\begin{equation*} %\label{eq:WLLNProof}
  \begin{split}
  \Prb \left(   \frac{ \sum\limits_{k=1}^n (\xi_k -
      \mathbb{E}\xi_k)}{b_n}> \epsilon \right)
   \leq& \frac{1}{\epsilon^2 b_n^2 }
  \E \left( \sum\limits_{k=1}^n (\xi_k -
      \mathbb{E}\xi_k)    \right)^2
    \leq  \frac{1}{\epsilon^2 b_n^2} \sum\limits_{k=1}^n
     \mathrm{Var} \xi_k\\
    \leq& \frac{1}{\epsilon^2 b_n^2} \sum\limits_{k=1}^m
    \mathrm{Var} \xi_k
    +\frac{1}{\epsilon^2}\sum\limits_{k=m}^n
   \frac{\mathrm{Var} \xi_k}{b_k^2}\\
    \leq& \frac{1}{\epsilon^2 b_n^2} \sum\limits_{k=1}^m
    \mathrm{Var} \xi_k
    +\frac{1}{\epsilon^2}\sum\limits_{k=m}^\infty
   \frac{\mathrm{Var} \xi_k}{b_k^2}.
 \end{split}
\end{equation*}
Since  $b_n \rightarrow \infty$,  \eqref{eq:stongConverLNN} follows.
\end{proof}

The following central limit theorem is a special case of a more general
result for martingales; see, for instance \cite[Theorem 5.5.4(II)]{LiptserShiryayevBookMartingales}.
\begin{Lem}[CLT for Stochastic Integrals]\label{thm:MartingaleCLT}
  Let $\mathcal{S} = (\Omega, \mathcal{F}, \Prb, \{\mathcal{F}_t\}_{t
    \geq 0},$ $\{W_{k}\}_{k \geq 1})$
  be a stochastic basis.  Suppose that
  $\sigma_k \in L^2(\Omega; L^2([0,T]))$ is a sequence of real valued
  predictable processes such that
  \begin{equation*}\label{eq:ConvScaleExactTo1}
    \lim_{N \rightarrow \infty} \frac{\sum\limits_{k =1}^{N} \int_0^T \sigma_k^2 dt}
      {\sum\limits_{k =1}^{N} \E \int_0^T \sigma_k^2 dt}   = 1
      \quad \textrm{ in Probability.}
  \end{equation*}
  Then
  \begin{equation*}%\label{eq:ConvCLTConclusion}
      \frac{\sum\limits_{k =1}^{N} \int_0^T \sigma_k dW_k}
      {\left(\sum\limits_{k =1}^{N} \E \int_0^T \sigma_k^2 dt\right)^{1/2}}
 \end{equation*}
 converges in distribution to a standard normal random variable
as $N \rightarrow \infty$.
\end{Lem}

\subsection{The Nonlinear Term}

The nonlinear term appearing in \eqref{eq:SNSEabs} is given by
$B(U,U^\sharp) := P_H ( (U \cdot \nabla) U^\sharp) = P_H (\sum_{j=1}^{2} U_j \pd{j} U^\sharp)$, which is
defined for $U \in V$ and $U^\sharp \in D(A)$. Note that, for brevity of
notation, we will often write $B(U)$ for $B(U,U)$ as for example in
\eqref{eq:SNSEabs}.
We have the following %classical
properties of $B$:

\begin{Lem}
  \mbox{}
  \begin{itemize}
  \item[(i)]   $B$ is bilinear and continuous from $V \times V$ into $V'$ and from
    $V \times D(A)$ into $H$.
    For $U, U^\sharp \in V$, $B$ satisfies the cancelation property
    \begin{equation}\label{eq:BCancel}
      \langle B(U, U^\sharp), U^{\sharp} \rangle = 0.
    \end{equation}
    If $U, U^\sharp, U^\flat$ are elements in $V$, then
    \begin{equation}\label{eq:BweakEstmates}
      \left| \langle B(U, U^{\sharp}), U^\flat  \rangle \right|
      \leq c |U|^{1/2}\| U \|^{1/2} \|U ^{\sharp}\| |U^{\flat}|^{1/2} \| U^{\flat}\|^{1/2}.
    \end{equation}
  On the other hand if $U \in V$, $U^\sharp \in D(A),$ and $U^\flat
  \in H,$ then we have:
  \begin{equation}\label{eq:BSizeEst}
      \left| (B(U,U^\sharp),U^\flat) \right|
        \leq c
        \begin{cases}
         |U|^{1/2}\|U\|^{1/2} \|U^\sharp\|^{1/2}
        |AU^\sharp|^{1/2}|U^\flat|.\\
         |U|^{1/2}|AU|^{1/2} \|U^\sharp\||U^\flat|.\\
      \end{cases}
  \end{equation}

  \item[(ii)] In the case of either periodic, (\ref{eq:PeriodicBC}) or
    Dirichlet, (\ref{eq:BCDirchlet}) boundary conditions, $B(U) \in D(A^{\beta})$
    for every $0< \beta < 1/4$ and every $U \in D(A)$.  Moreover, for such
    values of $\beta$,
    \begin{equation}\label{eq:SmallFracOrderEstimate}
      |A^{\beta} B(U)|^2 \leq c\|U\|^2 |AU|^2.
    \end{equation}
  \item[(iii)] In the case of periodic boundary conditions (\ref{eq:PeriodicBC}),
    $B(U,U^\sharp) \in D(A^{\beta})$ whenever $\beta > 1/2$,
    $U \in D(A^{\beta})$, $U^{\sharp} \in D(A^{\beta + 1/2})$, and for such
    $U$, $U^{\sharp}$,
    \begin{equation} \label{eq:HigherOrderEstNonlinearTerm}
      |A^{\beta} B( U, U^{\sharp})|^2 \leq c |A^\beta U|^2 |A^{\beta + 1/2} U|^2.
    \end{equation}
  \end{itemize}
\end{Lem}

\begin{proof}
  The properties outlined in (i) and (iii) are classical; see, for instance,  \cite{Temam1}, or
  \cite[Lemma 10.4]{ConstantinFoias1} for (\ref{eq:HigherOrderEstNonlinearTerm}).

  The properties in (ii) are established via  interpolation and the equivalence of certain
  fractional order spaces, see \cite{GlattHoltzTemam1}.  Since \cite{GlattHoltzTemam1}
  emphasized the case of spatial dimension $3$, for the sake completeness, we briefly
  recall the arguments.

  For any element $U \in D(A)$, standard estimates imply that
  \begin{equation*}%\label{eq:BEstL2H1}
    \begin{split}
      | B(U) |^2 &\leq c \|U\|^3 |A U| , \\
      \| B(U) \|_{H^1(\mathcal{M})^2}^2 &\leq c\|U\| |A U|^3.
    \end{split}
  \end{equation*}
  Let $\tilde{V} = H \cap H^1(\mathcal{D})^2$ and, for $s \in (0,1)$
  we define the interpolation spaces $\tilde{V}_{s} = [\tilde{V},H]_{1 - s}$.
  See \cite{LionsMagenes1} for the general theory.  In \cite{GlattHoltzTemam1},
  it is established that $D(A^\beta) = \tilde{V}_{2\beta}$, $\beta < 1/4$
  in Dirichlet case \eqref{eq:BCDirchlet}.\footnote{In the periodic case,
  \eqref{eq:PeriodicBC}, $\tilde{V} = V$ so that $D(A^\beta) = \tilde{V}_{2\beta}$,
  as a direct consequence of the fact that $D(A^{1/2}) = V$.}
  Note that $\tilde{V}$ does not incorporate
  boundary conditions and so $B(U) \in \tilde{V}$, for $U \in D(A)$. In consequence,
  for any such $U \in D(A)$ and allowed values of  $0 < \beta < 1/4$ we have, by interpolation
  \begin{equation*}%\label{eq:smallfracAEstB}
    \begin{split}
      |A^{\beta}B(U)|^{2} &= |B(U)|_{\tilde{V}^{2\beta}}^2
      \leq (|B(U)|^{1 - 2 \beta}\|B(U)\|^{2\beta}_{H^1})^2\\
      &\leq c(\|U\|^3 |A U|)^{1 - 2 \beta} (\|U\| |A U|^3)^{2\beta}\\
      &\leq c(\|U\|^3 |A U|)^{1/2} (\|U\| |A U|^3)^{1/2}
      \leq c\|U\|^2 |A U|^2.
    \end{split}
  \end{equation*}
  Combining these observations gives (ii), completing the proof.
\end{proof}

\begin{Rmk}
\label{rmk:ProblemsProblemsProblemsWithTheBoundary}
When we consider the case \eqref{eq:BCDirchlet}
it is not true in general that  $B(U) \in D(A^{\beta})$,
even for $U \in D(A^{\beta + 1/2})$, $\beta \geq 1/4$.
This is due to the fact that while the
Leray projector $P_{H}$ is continuous on $H^{m}(\mathcal{D})$,
$m \geq 1$,
we do not expect that $P_{H}$ maps $H^{m}_{0}(\mathcal{D})$
into $H^{m}_{0}(\mathcal{D})$.  See \cite{Temam1} and also
\cite{GlattHoltzTemam1}.  For this reason we may not expect
an inequality like \eqref{eq:HigherOrderEstNonlinearTerm}
for such Dirichlet boundary conditions.  As such,
\eqref{eq:SmallFracOrderEstimate} relies on a delicate
analysis of small fractional order space where the boundary
is not present;  see \cite{GlattHoltzTemam1, LionsMagenes1}.
\end{Rmk}

\subsection{Existence, Uniqueness and Higher Regularity}

With these mathematical formalities in place we now define precisely
(\ref{eq:SNSEabs}), in the usual time integrated sense and recall some
now well established existence, uniqueness and regularity results for
these equations. Note that for this work the solutions we consider correspond to
so called `strong solutions' in the deterministic setting (see \cite{Temam1}).
In the context of stochastic analysis, since we may suppose that
the stochastic basis $\mathcal{S}$ is fixed in advance, we may
say that the solutions considered are `strong' (or less confusingly
`pathwise') in the probabilistic sense as well.
\begin{Thm}\label{thm:WellPosedNessandHigherRegularity}
\mbox{}
\begin{itemize}
\item[(i)] Suppose that we impose either (\ref{eq:PeriodicBC}) or
  (\ref{eq:BCDirchlet}) and assume that $U_0 \in V$,
  $\sigma \in L_2(H, V)$.
  Then there exists a unique, $H$-valued, $\mathcal{F}_t$-adapted
  process $U$ with
  \begin{equation}\label{eq:StrongSolutionReg}
    U \in L^2_{loc} ([0,\infty);D(A)) \cap C([0,\infty);V)
    \quad a.s.
  \end{equation}
  and so that for each $t \geq 0$,
  \begin{equation*}%\label{eq:StrongSolutionEqn}
    U(t) + \int_0^t ( \nu AU + B(U) )dt'
    =U_0 + \sum_k \sigma \Phi_k W^k(t),
  \end{equation*}
  with the equality understood in $H$.
\item[(ii)] In the case of periodic boundary conditions
  (\ref{eq:PeriodicBC}) if $\beta > 1/2$ so that
  $\sigma \in L_2(H,D(A^\beta))$, $U_0 \in D(A^{\beta})$,
  then
  \begin{equation}\label{eq:HigherRegularityOnUPer}
    U \in L^2_{loc} ([0,\infty), D(A^{\beta+1/2})) \cap C([0,\infty), D(A^{\beta})).
  \end{equation}
\end{itemize}
\end{Thm}
\begin{Rmk}
\mbox{}
\begin{itemize}
\item[(i)] As noted above, when $\sigma$ is defined via
(\ref{eq:sigmaOperation}), $\sigma \in L_2(H, V)$ whenever $\gamma >1$.
Indeed we have $\sigma \in L_2(H, D(A^\beta))$ for every $\beta < \gamma-1/2$.
\item[(ii)] We suspect that higher regularity similar to Theorem \ref{thm:WellPosedNessandHigherRegularity},
(ii)  may be established
in the case of Dirichlet boundary conditions, \eqref{eq:BCDirchlet}.  However
since \eqref{eq:HigherOrderEstNonlinearTerm} does not apply
(see Remark~\ref{rmk:ProblemsProblemsProblemsWithTheBoundary})
a different proof than outlined here is needed.
\end{itemize}
\end{Rmk}
\begin{proof}
  The well-posedness of (\ref{eq:SNSEabs}) has been studied by many
  authors as discussed in the introduction.  Since we are considering
  the case of an additive noise the proof is close to the
  deterministic case after we perform a suitable change of variables.
  For completeness, we briefly recall some of the formal arguments and note
  that the computations may be rigorously justified with a
  suitable Galerkin scheme.  Consider first the linear system
  (\ref{eq:StokesStochasticAbs}) with initial condition $\bar{U}(0) = U_0$.
  As in Lemma~\ref{thm:RegStokesEqn} above, we have that
  $\bar{U}$ in $L^2_{loc} ([0,\infty);D(A)) \cap C([0,\infty);V)$
  (or in  $L^2_{loc} ([0,\infty),
  D(A^{\beta+1/2})) \cap C([0,\infty), D(A^{\beta}))$, under the
  conditions of item (ii)).
  We now consider the shifted variable $\tilde{U} = U - \bar{U}$, which satisfies
  \begin{equation}\label{eq:randPDEfromNoiseRemoval}
       \frac{d\tilde{U}}{dt} + \nu A\tilde{U} + B(\tilde{U} + \bar{U})
       = 0 \quad \tilde{U}(0) = 0\, .
  \end{equation}
  The estimates that lead to (\ref{eq:HigherRegularityOnUPer}) are
  standard. We first multiply (\ref{eq:randPDEfromNoiseRemoval}) by $U$,
  integrate over the domain and use
  (\ref{eq:BCancel}), (\ref{eq:BweakEstmates}), \eqref{eq:LinearEquationReg}
  to infer that $\tilde{U} \in L^2_{loc}([0,\infty); V) \cap L^\infty_{loc}([0,\infty);H)$.  With this regularity
  in hand we next
  multiply (\ref{eq:randPDEfromNoiseRemoval}) by $A \tilde{U}$
  and apply (\ref{eq:BSizeEst}),  \eqref{eq:LinearEquationReg}
  in order to conclude (\ref{eq:StrongSolutionReg}).

  For $\beta > 1/2$ we multiply (\ref{eq:randPDEfromNoiseRemoval}) by $A^{2\beta} \tilde{U}$and infer
 \begin{equation}\label{eq:HigherRegularityEstTilde}
    \frac{d |A^\beta \tilde{U}|^2}{dt}
    + 2\nu |A^{\beta +  1/2} \tilde{U}|^2  - 2\langle A^{\beta}
    B(\tilde{U} + \bar{U}), A^{\beta}U \rangle = 0.
 \end{equation}
 Since $\beta > 1/2$, we may apply
 (\ref{eq:HigherOrderEstNonlinearTerm}) and estimate
 \begin{equation*}%\label{eq:HigherRegEstYoung}
   \begin{split}
     \frac{d |A^\beta \tilde{U}|^2}{dt}
   +& 2\nu |A^{\beta +  1/2} \tilde{U}|^2 \\
   \leq& c |  A^\beta (\tilde{U}+ \bar{U})|
        |A^{\beta + 1/2}(\tilde{U} + \bar{U})|
        |A^\beta \tilde{U}|\\
   \leq& c (|  A^\beta \tilde{U}|^2+ |A^\beta\bar{U}|^2) |A^\beta \tilde{U}|^2
        +\nu |A^{\beta + 1/2}\tilde{U}|^2 + \nu |A^{\beta + 1/2}\bar{U}|^2 .
  \end{split}
\end{equation*}
 Rearranging,
 \begin{equation*}%\label{eq:RegularityEstConclusion}
       \frac{d |A^\beta \tilde{U}|^2}{dt}
   + 2\nu |A^{\beta +  1/2} \tilde{U}|^2 \leq
   c (|  A^\beta \tilde{U}|^2+ |A^\beta\bar{U}|^2) |A^\beta \tilde{U}|^2
        + \nu |A^{\beta + 1/2}\bar{U}|^2 .
 \end{equation*}
 Observe that, due to the Gronwall Lemma, if $\tilde{U} \in
 L^2_{loc}([0, \infty); D(A^{\beta}))$
 then we infer that $\tilde{U}
 \in L^2_{loc}([0, \infty); D(A^{\beta+1/2})) \cap L^\infty_{loc}([0,
 \infty); D(A^{\beta}))$.
 The desired result therefore follows from an inductive argument on $\beta$
 starting with the base case assumption $\beta \in [1/2, 1)$ which
 is satisfied as a consequence of \eqref{eq:StrongSolutionReg}.
\end{proof}

\section{Estimators for $\nu$: Heuristic Derivation and the Main Results}
\label{sec:DerivationOfTheEstimator}
In this section we sketch the heuristic derivations of the estimators
based on a particular version of the Girsanov Theorem.
We then restate, now in precise terms, the main
results of this paper.

As before we denote by $U^N$ the projection of the solution $U$ of the original equation \eqref{eq:SNSEabs} onto
$H_N = P_N H \cong \mathbb{R}^{N}$.
Note that $U^N$ satisfies the following  finite dimensional system:
\begin{equation}\label{eq:SNSEabsFinite}
  dU^N =  -(\nu AU^N + \psi^N)dt +  P_N\sigma dW,
  \quad U^N(0) = U_0^N,
\end{equation}
where $\psi^N(t):=P_N(B(U))$.
To obtain an initial guess of the form of the estimator for the parameter $\nu$, we
treat $\psi^N$ as an external known quantity, independent of $\nu$
and view \eqref{eq:SNSEabsFinite}
as a stochastic equation evolving in $\mathbb{R}^N$.
Let us denote by $\mathbb{P}^{N,T}_\nu$ the probability measure in $C([0,T]; \mathbb{R}^N)$
generated by $U^N$.  Formally, we compute the Radon-Nikodym derivative or
Likelihood Ratio $d\Prb_{\nu}^{N,T}/d\Prb_{\nu_0}^{N,T}$ (see e.g. \cite[Section 7.6.4]{LiptserShiryayev})
\begin{equation*}
  \begin{split}
    \frac{d\mathbb{P}^{N,T}_{\nu} (U^N) } {d\mathbb{P}^{N,T}_{\nu_0}}
    =   \exp\Big( & -\int_0^T ( \nu -\nu_0)(AU^N )^\prime G^2 dU^N(t) \\
       & -  \frac{1}{2}\int_0^T (\nu^2 -\nu_0^2) (AU^N)^\prime G^2 AU^N dt \\
       & - \int_0^T (\nu-\nu_0)(AU^N)^\prime G^2 \psi^Ndt \Big),
 \end{split}
\end{equation*}
where $G := (P_N \sigma)^{-1} = \mathrm{diag}[\sigma_1^{-1}, \ldots, \sigma_N^{-1}] =
\mathrm{diag}[\lambda_1^{\gamma}, \ldots, \lambda_N^{\gamma}] $
and $v^\prime$ denotes the transpose of the vector $v\in\mathbb{R}^N$.
By maximizing the Likelihood Ratio with respect to the parameter of interest
$\nu$, we may compute the (formal) Maximum Likelihood Estimator (MLE)
$\nu_N$ of the parameter $\nu$.
A direct computation yields
\begin{equation}\label{eq:formalMLEVectorLang}
  \nu_{N} =
  -\frac{\int_0^T (AU^N)^\prime G^2 dU^N  + \int_0^T (AU^N)^\prime G^2 P_N(B(U)) dt  }
  { \int_0^T (AU^N)^\prime G^2 AU^N dt}.
\end{equation}
As expected, $\nu_{N}$ is a valid estimator and in fact one can
show that it is a consistent estimator of the true parameter $\nu$.
This consistency makes essential use of the fact that
the denominator  $\int_0^T (AU^N)^\prime G^2 AU^N dt$ diverges to infinity as $N \uparrow \infty$.
See Lemma~\ref{thm:ExactEst1} below.
With this in mind, we introduce a slight modification to the MLE \eqref{eq:formalMLEVectorLang}, and
propose  the following class of estimators
\begin{equation}\label{eq:EstimatorNuTilde}% \label{eq:EstimatorNuTilde}
 \widetilde{\nu}_N =
 -\frac{\int_0^T \langle A^{1 + 2 \alpha} U^N, dU^N \rangle
   + \int_0^T \langle A^{1 + 2 \alpha}U^N, P_N B(U) \rangle dt  }
 { \int_0^T |A^{1+ \alpha}U^N|^2 dt},
\end{equation}
where $\alpha$ is a free parameter with a range specified later on.
Note that this formulation appears in the functional language developed above
and is derived using that
the action of $G^2$ on $H_N$ is equivalent to $A^{2\gamma}$.
Also we observe that $\nu_{N}$ is a particular case of $\widetilde{\nu}_N$ with
$\alpha=\gamma$.

While the estimator $\widetilde{\nu}_N$ has desirable
theoretical properties, it also assumes that $P_N(B(U))$ is computable,
which could be quiet a difficult task.
Since our goal is to provide estimators that can be eventually implemented in practice
(evaluated numerically),
we propose two further classes of estimators. One class is naturally derived from
\eqref{eq:EstimatorNuTilde} by
approximating $P_N(B(U))$ with $P_N(B(U^N))$
\begin{equation}\label{eq:EstimatorNuCheck}
\check{\nu}_N = -  \frac{\int_0^T \langle A^{1 + 2 \alpha} U^N, dU^N \rangle
   + \int_0^T \langle A^{1 + 2 \alpha}U^N, P_N B(U^N) \rangle dt  }
 { \int_0^T |A^{1+ \alpha}U^N|^2 dt}.
\end{equation}
Note that $\check{\nu}_N$ now depends only on the first $N$ Fourier modes.
However, even in this case the
expression for $P_N B(U^N)$ is very complicated due to the
nontrivial coupling of the modes.  See e.g. \cite{DoeringGibbon1}.
It turns out, as shown rigorously below
(see Proposition~\ref{thm:ConsistencyInNonlinearTermPlusRates}),
that the second term appearing in \eqref{eq:EstimatorNuTilde},
\begin{equation}\label{eq:ErrorTerm}
  \kappa_N := - \frac{ \int_0^T \langle A^{1 + 2 \alpha}U^N, P_N B(U) \rangle dt  }
 { \int_0^T |A^{1+ \alpha}U^N|^2 dt},
\end{equation}
is of lower order and tends to zero, as $N\to\infty$.
Hence we get the following consistent estimators of the parameter $\nu$
\begin{equation}\label{eq:EstimatorNuHat}
  \begin{split}
  \hat{\nu}_N
   &= -\frac{\int_0^T \langle A^{1 + 2\alpha} U^N, d U^N \rangle}
         {\int_0^T |A^{1+ \alpha } U^N|^2dt}
    = -
    \frac{\sum\limits_{k =1}^{N} \lambda_k^{1 + 2 \alpha} \int_0^T u_k du_k}
        {\sum\limits_{k =1}^{N}\lambda_k^{2 + 2 \alpha} \int_0^T u_k^2 dt}\\
   &= - \frac{\sum\limits_{k =1}^{N} \lambda_k^{1 + 2 \alpha}
           (u_k^2(T) - u_k^2(0) - T \lambda_k^{-2\gamma}) }
        {2 \sum\limits_{k =1}^{N} \lambda_k^{2 + 2 \alpha} \int_0^T u_k^2 dt}.\\
  \end{split}
\end{equation}
Clearly this last estimator is easiest to compute numerically.  On the other hand
it may lack the speed of convergence of the first two.

We conclude this section with
the main result of this paper:
\begin{Thm}\label{thm:MainResultTechincal}
  Suppose that $U$ solves \eqref{eq:SNSEabs} with either (\ref{eq:PeriodicBC}) or (\ref{eq:BCDirchlet})
  in the sense of and under the conditions imposed by
  Theorem~\ref{thm:WellPosedNessandHigherRegularity}.  Assume that
  $\gamma > 1$ and in the case  (\ref{eq:BCDirchlet}), additionally that
  $\gamma < 1 + 1/4$.   Also, assume that $U_0 \in D(A^\beta)$, for some $\beta > \gamma - 1/2$.
  \begin{itemize}
  \item[(i)] If $\alpha>\gamma-1$, then $\widetilde{\nu}_N, \ \check{\nu}_N$ and $\hat{\nu}_N$
  as given by \eqref{eq:EstimatorNuTilde}, \eqref{eq:EstimatorNuCheck},  and
  \eqref{eq:EstimatorNuHat}
  are weakly consistent estimators of the parameter $\nu$, i.e.
    \begin{equation*}%\label{eq:consistencynuHat}
      \lim_{N \rightarrow \infty} \widetilde{\nu}_N
      = \lim_{N \rightarrow \infty} \check{\nu}_N
      = \lim_{N \rightarrow \infty} \hat{\nu}_N  = \nu
   \end{equation*}
    in probability.
  \item[(ii)] If  $\alpha > \gamma-1/2$, then
    $\tilde{\nu}_N$ is asymptotically normal with rate $N$, i.e.
    \begin{equation}\label{eq:asymtoticNormalnuHat}
     N (\tilde{\nu}_N -\nu)
        \overset{d}{\longrightarrow} \eta\, ,
  \end{equation}
where $\eta$ is Gaussian random variable with mean zero and variance\\
 $\frac{2\nu(\alpha-\gamma+1)^2}{\lambda_1T(\alpha-\gamma+1/2)}$.
  \end{itemize}
\end{Thm}

\section{Proof of the Main Theorem}
\label{sec:ProofOfMainRes}

We establish the proof of the Theorem~\ref{thm:MainResultBasic} in
a series of propositions.
As mentioned in the introduction, we do not have precise
spectral information about Fourier coefficients $u_k = (U, \Phi_k), \ k\geq 1$, in contrast to linear case (see Section \ref{sec:StokesEqAndLimitThrms}).
To overcome this, we proceed by decomposing the solution into
a linear and a nonlinear part,
$U = \bar{U} + R$.  We assume that $\bar{U}$ is the solution of
the linear stochastic Stokes equation \eqref{eq:StokesStochasticAbs}
with $\bar{U}(0) = 0$.  The residual $R$ must therefore satisfy,
\begin{equation}\label{eq:resEqn}
 \pd{t} R + \nu A R  = - B(U),  \quad R(0) = R_0.
\end{equation}

First, we study the regularity properties
of $R$ and show that $R$ is slightly smoother than $\bar{U}$.
Subsequently, we make crucial use of this extra regularity and
establish the consistency of the proposed estimators by showing
that second term in \eqref{eq:EstimatorNuTilde} converges to zero.
The final section treats the asymptotic
normality using CLT introduced
in Section~\ref{sec:StokesEqAndLimitThrms}.

\begin{Rmk}\label{rmk:PresentationClairy}
   For simplicity and clarity of presentation we shall
   assume a more regular initial condition $U_{0} \in D(A^{\gamma})$
   in contrast to the statement of Theorem~\ref{thm:MainResultTechincal}.
   The more general case when we assume merely that $U_{0} \in D(A^{\beta})$
   for some $\beta > \gamma -1/2$ may be treated by writing $U = \bar{U} + R
   +S$, where $\bar{U}$ satisfies \eqref{eq:StokesStochasticAbs}
   with $\bar{U}_{0} = 0$, $R$ satisfies \eqref{eq:resEqn}, this
   time with $R_{0} = 0$, and finally $S$ is the solution of
   $\pd{t} S + \nu A S = 0$, with $S(0) = U_{0}$.
\end{Rmk}

\subsection{Regularity Properties for the Residual}
\label{sec:residual-estimates}

\begin{Prop}\label{thm:RegForRes}
  Suppose that $\beta\geq 0$, $\gamma > 1$
  and that $R$ solves \eqref{eq:resEqn} with $U$
  the solution of \eqref{eq:SNSEabs} corresponding
  to an initial condition $U_0 \in D(A^{1/2 + \beta})$.
  \begin{itemize}
  \item[(i)] If $U$ and $R$ satisfy Dirichlet boundary conditions \eqref{eq:BCDirchlet}, and
  $\beta < 1/4$, then for every $T > 0$ we have
    \begin{equation}\label{eq:ptwiseHigherRegrho}
      \sup_{t \in [0,T]} | A^{1/2 + \beta} R |^{2} +
      \int_{0}^{T} |A^{1 + \beta} R |^{2} < \infty.
    \end{equation}
    Moreover, for an increasing sequence of stopping times $\tau_n$ with
    $\tau_n \uparrow \infty$,
    \begin{equation}\label{eq:StoppedEValuEst1}
      \E \left( \sup_{t \in [0,\tau_{n}]} | A^{1/2 + \beta} R |^{2} +
      \int_{0}^{\tau_n} |A^{1 + \beta} R |^{2} \right)< \infty.
    \end{equation}
  \item[(ii)] In the case that both $U$ and $R$
    satisfy periodic boundary conditions \eqref{eq:PeriodicBC}
  and we assume $\beta < \gamma - 1/2$ the same
    conclusions hold.
  \end{itemize}
\end{Prop}

\begin{proof}
As above in Theorem~\ref{thm:WellPosedNessandHigherRegularity}
the computations given here may be rigoursly justified via Galerkin
approximations.  Multiplying \eqref{eq:resEqn} by $A^{1 + 2\beta} R$,
integrating and using the symmetry of the powers of $A$ we infer
\begin{equation}\label{eq:resEqnEngry}
 \frac{1}{2} \frac{d}{dt} |A^{\beta + 1/2} R |^{2} + \nu |A^{\beta+1} R
 |^{2}  = - \langle A^{\beta}B(U), A^{\beta+ 1} R  \rangle.
\end{equation}
For the case of a bounded domain, \eqref{eq:BCDirchlet}, we infer
from \eqref{eq:SmallFracOrderEstimate} and
Theorem~\ref{thm:WellPosedNessandHigherRegularity},
\eqref{eq:StrongSolutionReg} that,
\begin{equation*} % \label{eq:uRegularity}
  \int_{0}^{T}|A^{\beta}B(U)|^{2} dt  \leq c \int_0^T \|U\|^2 |AU|^2 dt < \infty  \quad a.s.
\end{equation*}
By integrating \eqref{eq:resEqnEngry} in time and making standard estimates
with Young's Inequality, \eqref{eq:ptwiseHigherRegrho} now follows
in this case.

In the case of the periodic domain we estimate $|A^{\beta}B(U)|$ differently.
Define $\beta' = \max\{\beta, \gamma/2\}$ so
that $1/2 < \beta' < \gamma - 1/2$.
By applying the higher regularity
estimates (\ref{eq:HigherOrderEstNonlinearTerm}) we find that
\begin{equation*}%\label{eq:BestforResPeriodicDomain}
  \begin{split}
    | \langle A^{\beta} B(U), A^{\beta + 1}R \rangle |
     \leq& |A^{\beta} B(U)||A^{\beta +1}R |\\
     \leq& c |A^{\beta'} B(U)||A^{\beta +1}R |\\
     \leq& c |A^{\beta'} U||A^{\beta' +1/2} U| |A^{\beta +1}R |\\
     \leq& c |A^{\beta'} U|^2|A^{\beta' +1/2} U|^2
                 + \frac{\nu}{2}|A^{\beta +1}R |^2 \, .
  \end{split}
\end{equation*}
Due to Theorem~\ref{thm:WellPosedNessandHigherRegularity},
(ii), we have, for any $T >0$, that
\begin{equation*}%\label{eq:RegTranlation}
  \int_0^T |A^{\beta'} U|^2|A^{\beta' +1/2} U|^2 dt < \infty
  \quad a.s.,
\end{equation*}
and \eqref{eq:ptwiseHigherRegrho} follows once again.

For the stopping times $\tau_n$, we define
\begin{equation*}\label{eq:stoppingTimesBndDomain}
  \tau_n := \inf_{t \geq 0}
  \left\{ \sup_{t' \leq t} \|U\|^2
    + \int_0^t |AU|^2dt' > n
   \right\}
\end{equation*}
when \eqref{eq:BCDirchlet} is assumed and
\begin{equation*}\label{eq:stoppingTimesBndDomain}
  \tau_n := \inf_{t \geq 0}
  \left\{ \sup_{t' \leq t} |A^{\beta'} U|^2
    + \int_0^t |A^{\beta' + 1/2} U|^2dt' > n
   \right\}
\end{equation*}
for \eqref{eq:PeriodicBC}.  In either case it is clear that
$\{\tau_{n}\}_{n \geq 1}$ is increasing.  Moreover, in
 the case \eqref{eq:BCDirchlet}, since $\Prb( \tau_{n}  < T)
 = \Prb( \sup_{t' \leq T} \|U\|^2
    + \int_0^T |AU|^2dt' \geq n)$, it follows from
   \eqref{eq:StrongSolutionReg} and the fact that
   $\tau_{n}$ is increasing, that $\lim_{n \rightarrow}\tau_{n}=\infty$
   a.s.  Arguing in the same manner for the case
   \eqref{eq:PeriodicBC}, the proof is complete.
\end{proof}

\begin{Rmk*}\label{rmk:MoreOnRegularity}
  Comparing Proposition~\ref{thm:RegForRes},
  \eqref{eq:ptwiseHigherRegrho} with Lemma~\ref{thm:RegStokesEqn},
  \eqref{eq:LinearEquationReg}, \eqref{eq:BlowUpRateLinearSystem}
  we see that $R$ has been shown to be just shy of a derivative
  more regular than $\bar{U}$.  More precisely we have that,
  for any $\epsilon >0$,
  \begin{displaymath}
    \bar{U} \in
      L^{2}(\Omega; L^{2}_{loc}([0,\infty); D(A^{\gamma- \epsilon}))),
      \quad
          \bar{U} \not \in
      L^{2}(\Omega; L^{2}_{loc}([0,\infty); D(A^{\gamma + \epsilon}))),
  \end{displaymath}
  while on the other hand,
  \begin{displaymath}
      R \in
      L^{2}(\Omega; L^{2}_{loc}([0,\infty); D(A^{\gamma + 1/2 - \epsilon}))).
  \end{displaymath}
\end{Rmk*}

As an immediate application of these properties of the residual $R$
we have the following result:
\begin{Lem}\label{thm:ExactEst1}
  Suppose that $U$ and $\bar{U}$ are the solutions
  of (\ref{eq:SNSEabs}) and (\ref{eq:StokesStochasticAbs})
  respectively.  For both  \eqref{eq:PeriodicBC},  and \eqref{eq:BCDirchlet}
  we suppose that $\gamma > 1$, $U(0) = U_0 \in D(A^{\gamma})$
        \footnote{At the cost of further evaluations,
        this condition may be weakened to the conditions imposed
        in Theorem~\ref{thm:MainResultTechincal}.  This applies both here and below for Propositions
        \ref{thm:ConsistenccyIntheStochasticTerms}, \ref{thm:ConsistencyInNonlinearTermPlusRates}.
        See Remark~\ref{rmk:PresentationClairy}.}
  and $\bar{U}_0 = 0$.  Additionally,
  in the case (\ref{eq:BCDirchlet}), we assume that $\gamma <1+1/4$.
  Then, for any $\alpha > \gamma - 1$,
 \begin{equation} \label{eq:ExactEst1}
    \lim_{N \rightarrow \infty}
    \frac {\int_0^T  | A^{1 + \alpha } U^N|^{2} dt}
        { \E \int_0^T  | A^{1 + \alpha } \bar{U}^N|^{2} dt} = 1
   \end{equation}
with probability one.
\end{Lem}
\begin{proof} Note that
\begin{align*}
|A^{1+ \alpha} U^N|^2
     \leq& |A^{1+ \alpha} \bar{U}^N|^2 + |A^{1+ \alpha} R^N|^2
    + 2 |A^{1+ \alpha} \bar{U}^N||A^{1+ \alpha} R^N|,\\
|A^{1+ \alpha} U^N|^2
     \geq & |A^{1+ \alpha} \bar{U}^N|^2 + |A^{1+ \alpha} R^N|^2
    - 2 |A^{1+ \alpha} \bar{U}^N||A^{1+ \alpha} R^N|\, ,
\end{align*}
and therefore \eqref{eq:ExactEst1} follows once we have shown that
\begin{equation*}%\label{eq:I1NConvLLN}
   I_{1}^{N} :=  \frac { \int_0^T  | A^{1 + \alpha } \bar{U}^N|^{2}}{  { \E \int_0^T  | A^{1 + \alpha } \bar{U}^N|^{2} } }
   \rightarrow 1 \quad a.s.
\end{equation*}
and that
\begin{equation}\label{eq:I2NConvReg}
   I_{2}^{N} :=  \frac { \int_0^T | A^{1 + \alpha }R^{N}|^{2} dt }{  { \E \int_0^T  | A^{1 + \alpha } \bar{U}^N|^{2} }}
   \rightarrow 0 \quad a.s.
\end{equation}
For the first item, $I_{1}^{N}$, we apply the law of large number (LLN),
Lemma~\ref{thm:LLN}, with $\xi_n:=\lambda_n^{2\alpha+2} \int_0^T \bar{u}_n^2(t)dt$
and $b_n := \sum\limits_{k=1}^n\mathbb{E}[\xi_k]$. Notice that, due to, \eqref{eq:LambdaOrder} and
(\ref{eq:specL2ubark})
\begin{equation}\label{eq:bnfliesoffthehandle}
  b_n
 \sim \sum_{k=1}^n \lambda_k^{2\alpha+2} \lambda_k^{-1 - 2\gamma}
 \sim \sum_{k =1}^n k^{2\alpha- 2 \gamma +1}.
\end{equation}
Given the assumptions  $\alpha>\gamma-1$, we have that
$\lim_{n \rightarrow \infty} b_n = \infty$.
Moreover, combining \eqref{eq:bnfliesoffthehandle} with
\eqref{eq:LambdaOrder},
\eqref{eq:specL4Typeubark}, (\ref{eq:partialSum})
\begin{equation*} %  \label{eq:LLNCalcVarL4ahh}
  \begin{split}
    \sum\limits_{n\geq 1} \frac{\mathrm{Var}\xi_n  }{b_n^2}
     &\sim \sum\limits_{n\geq 1}
        \frac{\lambda_n^{4\alpha + 4 - 3 - 4\gamma}}
       { \Big( \sum\limits_{k=1}^n \lambda_k^{2\alpha-2\gamma + 1} \Big)^2 } \\
   &\sim\sum\limits_{n\geq 1}
        \frac{\lambda_n^{4\alpha - 4\gamma + 1}}
        { (\lambda_n^{ 2\alpha-2\gamma + 2})^2} \sim \sum\limits_{n\geq 1} \frac{1}{ n^{3}} <\infty.
 \end{split}
\end{equation*}
Thus, by the LLN we conclude that $\lim_{N\to\infty} I_1^N=1$ with probability one.

Since $1 + \alpha > \gamma$, by \eqref{eq:BlowUpRateLinearSystem} we infer
\begin{equation*} % \label{eq:OnderOfDenominatorGrowth}
  \E\int_0^T |A^{1+\alpha} \bar{U}^N|dt \sim N^{2 \alpha - 2\gamma + 2} .
\end{equation*}
Pick any $\alpha' \in (\gamma - 1, \min\{\alpha, 1/4\})$, in the
case (\ref{eq:BCDirchlet}), or any $\alpha' \in (\gamma - 1,
\min\{\alpha, \gamma - 1/2\})$ under the assumption
(\ref{eq:PeriodicBC}).
By applying (\ref{eq:ptwiseHigherRegrho})
for $R$ established in Proposition~\ref{thm:RegForRes}, we have
in both cases that
\begin{equation*}%\label{eq:regConclusionAlphaPrime}
  \int_0^T | A^{1 + \alpha' } R|^{2} dt  < \infty  \quad \textrm{a.s.}
\end{equation*}
Combining these observations and making use of (\ref{eq:decompEstimates}), we have
\begin{align*}%\label{eq:squeezingOutRterms}
  I_{2}^{N} \leq& c
  \frac { \int_0^T | A^{1 + \alpha }R^{N}|^{2} dt }
      { N^{2\alpha -2 \gamma+ 2} }
     \leq  c
  \frac { \lambda_N^{2(\alpha - \alpha')} \int_0^T | A^{1 + \alpha' } R^{N}|^{2} dt }
     { N^{2\alpha -2 \gamma+ 2} }
      \leq  c
      \frac { \int_0^T | A^{1 + \alpha' } R|^{2} dt }
     { N^{2\alpha' -2 \gamma+ 2} } \, .
\end{align*}
Due to the restrictions on the choice of $\alpha'$, we have that $2
\alpha' - 2\gamma + 2 > 0$, and hence $I_2^N\to 0$, as $N\to\infty$, with probability one.
The proof is complete.
\end{proof}

\subsection{Consistency of the  Estimators}
\label{sec:cons-estim}

Using the dynamics of $U^N$, i.e. substituting \eqref{eq:SNSEabsFinite}
into \eqref{eq:EstimatorNuTilde}, we
get the following representation for the estimator $\widetilde{\nu}_N$:
\begin{equation}\label{eq:RepresentationNuTilde}
\begin{split}
  \widetilde{\nu}_N & =
  \nu  -
        \frac{\int_0^T \langle A^{1 + 2\alpha} U^N, P_N \sigma dW\rangle}
        {\int_0^T |A^{1+ \alpha } U^N|^2dt} \\
  & = \nu     -
        \frac{\int_0^T \langle A^{1 + 2\alpha - \gamma} U^N
          , \sum_{k=1}^N \Phi_k dW_k\rangle}
        {\int_0^T |A^{1+ \alpha } U^N|^2dt}.
\end{split}
\end{equation}
Similarly, we deduce
\begin{equation}\label{eq:RepresentationNuCheck}
\begin{split}
\check{\nu}_N=  \nu &-
        \frac{\int_0^T \langle A^{1 + 2\alpha - \gamma} U^N
          , \sum_{k=1}^N \Phi_k dW_k\rangle}
        {\int_0^T |A^{1+ \alpha } U^N|^2dt} \\
        & + \frac{ \int_0^T \langle A^{1 + 2 \alpha}U^N, P_N B(U)-P_NB(U^N) \rangle dt  }
 { \int_0^T |A^{1+ \alpha}U^N|^2 dt} .
 \end{split}
\end{equation}
Note that  $\hat{\nu}_N = \widetilde{\nu}_N -\kappa_N$, with $\kappa_N$ defined by
\eqref{eq:ErrorTerm}. Thus,
\begin{equation}\label{eq:RepresentationNuHat}
\widehat{\nu}_N = \nu - \kappa_N - \frac{\int_0^T \langle A^{1 + 2\alpha - \gamma} U^N
          , \sum_{k=1}^N \Phi_k dW_k\rangle}
        {\int_0^T |A^{1+ \alpha } U^N|^2dt}.
\end{equation}

With the above representations for the estimators, the consistency will follow if we show that
each stochastic term on the right-hand side of  \eqref{eq:RepresentationNuTilde},
\eqref{eq:RepresentationNuCheck}, \eqref{eq:RepresentationNuHat} converges to zero.

\begin{Prop}\label{thm:ConsistenccyIntheStochasticTerms}
  Assume the conditions and notations from Lemma~\ref{thm:ExactEst1}.
  Then,
  \begin{enumerate}
\item[(i)]  for every $\delta_{1} < \min \{2 + 2\alpha - 2 \gamma, 1 \}$,
 \begin{equation}\label{eq:effStochasticTermsLin}
 \lim_{N \rightarrow \infty} N^{\delta_{1}} \frac{\int_0^T \langle A^{1 + 2 \alpha -
    \gamma} \bar{U}^N, \sum_{k=1}^{N}  \Phi_{k} dW_{k} \rangle}
   { \int_0^T  | A^{1 + \alpha } U^N |^2 dt } = 0
   \quad \textrm{a.s.}
 \end{equation}
\item[(ii)]  whenever $\delta_{2} < \min \{2 + 2\alpha - 2 \gamma, 3/2 \}$ in
  the case \eqref{eq:PeriodicBC}, or whenever
  $\delta_{2} < \min \{2 + 2\alpha - 2 \gamma, 5/4  + 1 - \gamma \}$
  in the case \eqref{eq:BCDirchlet}, we have
  \begin{equation}\label{eq:effStochasticTermsResid}
    \lim_{N \rightarrow \infty} N^{\delta_{2}} \frac{\int_0^T \langle A^{1 + 2 \alpha -
     \gamma} R^N, \sum_{k=1}^{N}  \Phi_{k} dW_{k} \rangle}
    { \int_0^T  | A^{1 + \alpha } U^N |^2 dt } = 0
 \end{equation}
 in probability.
\end{enumerate}
\end{Prop}
\begin{proof}
  Due to Lemma~\ref{thm:ExactEst1}, (\ref{eq:ExactEst1}) and
  (\ref{eq:BlowUpRateLinearSystem}) the desired result follows once we
  show that each of sequences
  \begin{equation*}%\label{eq:J1Decomp}
      J_N^1 :=
      \frac{\int_0^T \langle A^{1 + 2 \alpha -
      \gamma}\bar{U}^N, \sum_{k=1}^{N}  \Phi_{k} dW_{k} \rangle}
      {\lambda_N^{2 + 2\alpha - 2 \gamma  - \delta_{1}}}
      =\frac{\sum_{k=1}^N \lambda^{1+2\alpha - \gamma}_k \int_0^T\bar{u}_k dW_k}
      {\lambda_N^{2 + 2\alpha - 2 \gamma  - \delta_{1}}}
  \end{equation*}
 and
 \begin{equation*} %  \label{eq:J2Decomp}
       J^{2}_{N} := \frac{\int_0^T \langle A^{1 + 2 \alpha -\gamma}
    R^{N}, \sum_{k=1}^{N}  \Phi_{k} dW_{k} \rangle}
  { \lambda_N^{2 + 2\alpha - 2 \gamma -\delta_{2}}}
 \end{equation*}
 converge to zero as $N \rightarrow \infty$.

 For the first term, $J^1_N$, define
 $\bar{\xi}_k := \lambda_k^{1 + 2\alpha -    \gamma} \int_0^T\bar{u}_k dW_k$ and
 $b_n := \lambda_n^{2 + 2\alpha - 2\gamma - \delta_{1}}$.  Under the given
 conditions, $\lim_{n \rightarrow \infty} b_n =\infty$.
 With the It\={o} Isometry and (\ref{eq:specL2ubark}), we have
 \begin{equation*} %\label{eq:VarukIndividual}
   \mathrm{Var} [\bar{\xi}_k] = \E [\bar{\xi}^2_k]
   \sim \lambda_k^{2 + 4 \alpha - 2\gamma}\lambda_k^{-(1+ 2\gamma)}
   =\lambda_k^{1 + 4 \alpha - 4\gamma}.
 \end{equation*}
Thus,
 \begin{equation*} % \label{eq:LLNApp1}
   \sum\limits_{n\geq 1}\frac{\mathrm{Var}\xi_n  }{b_n^2}
   \sim \sum\limits_{n\geq 1} \frac{\lambda_n^{1 + 4 \alpha - 4\gamma}}
          {\lambda_n^{4 + 4\alpha - 4 \gamma  - 2\delta_{1}}}
          = \sum\limits_{n\geq 1} \frac{1}{\lambda_n^{3 - 2 \delta_{1}}}
          \sim \sum\limits_{n\geq 1} \frac{1}{n^{3 - 2\delta_{1}}} <\infty.
  \end{equation*}
  Note that under the given conditions $\delta_{1} < 1$. This justifies
  the assertion that the final sum is finite.  We
  conclude, by the LLN, Lemma~\ref{thm:LLN} that
  $\lim_{N \rightarrow \infty} J_N^1 = 0$.

We turn to $J_N^2$.  Let $r_k := (R, \Phi_k), \ k\geq 1$, and for any stopping time $\tau$ we define
   \begin{equation*}%\label{eq:ResLLNRVDef}
      \zeta_k^\tau
      := \lambda_k^{1 + 2\alpha -\gamma} \int_0^\tau r_k dW_k .
  \end{equation*}
 Note that the random variables $\zeta_k^\tau, \ k\geq 1$, are uncorrelated.
 Similarly to the above arguments, we let $b_n := \lambda_n^{2 + 2\alpha - 2\gamma -
  \delta_{2}}$ and observe that this sequence is increasing and unbounded.
 Up to any stopping time $\tau$ such that $\mathrm{Var} \zeta_k^\tau < \infty$, we have
 \begin{equation}\label{eq:PrancarelType}
   \begin{split}
    \sum_{k \geq 1} \frac{\mathrm{Var} [\zeta_k^\tau]}{ b_k^2}
   &=  \sum_{k \geq 1}
        \frac{ \lambda_k^{2 + 4\alpha -2\gamma}}
        {\lambda_k^{4 +    4\alpha - 4 \gamma -2\delta_{2}}}
        \E \int_0^\tau r_k^2 dt\\
   &=  \sum_{k \geq 1}
        \lambda^{2\gamma - 2 + 2 \delta_{2}}_k
       \E \int_0^\tau r_k^2 dt\\
   &= \E \int_0^\tau | A^{\gamma  - 1 +  \delta_{2}} R|^2 dt \, .
   \end{split}
  \end{equation}
Note that under initial assumptions,
in the case of a bounded domain, \eqref{eq:BCDirchlet},
  $\gamma - 1 + \delta_{2} < 5/4$, and
  in the periodic case, \eqref{eq:PeriodicBC},
  we have $\gamma - 1 + \delta_{2} < \gamma + 1/2$.
  In either case, by taking $\tau_n$ as in
  Proposition~\ref{thm:RegForRes},  we infer from
  \eqref{eq:PrancarelType} with \eqref{eq:StoppedEValuEst1} that,
  for every $n$,
  \begin{equation*} % \label{eq:FiniteMeanTypeCondRealizedUpToStoppingTimes}
     \sum_{k \geq 1} \frac{\mathrm{Var} \zeta_k^{T \wedge \tau_n}}{ b_k^2}
       < \infty.
  \end{equation*}
  By applying Lemma~\ref{thm:LLN}, we conclude that, for each $n$ fixed,
  \begin{equation*} %\label{eq:StochasticIntegralResLLNConclusionST}
   \lim_{N \rightarrow \infty}
   \frac{\int_0^{T \wedge \tau_n} \langle A^{1 + 2 \alpha -\gamma}
     R^{N}, \sum_{k=1}^{N}  \Phi_{k} dW_{k} \rangle}
   { \lambda_N^{2 + 2\alpha - 2 \gamma -\delta_{2}}} = 0
   \quad \textrm{ in Probability.}
 \end{equation*}
 Since $\tau_n$ is  increasing, $\tilde{\Omega} =
  \cup_n \{ \tau_n > T \}$ is a set of full measure, and
  a simple estimate yields that
  \begin{equation*} %\label{eq:StochasticIntegralResLLNConclusion}
   \lim_{N \rightarrow \infty}
   \frac{\int_0^{T} \langle A^{1 + 2 \alpha -\gamma}
     R^{N}, \sum_{k=1}^{N}  \Phi_{k} dW_{k} \rangle}
   { \lambda_N^{2 + 2\alpha - 2 \gamma -\delta_{2}}} = 0
   \quad \textrm{ in Probability.}
 \end{equation*}
The proof is complete.
 \end{proof}

\begin{Cor}
Putting the admissible values $\delta_{1} = \delta_{2} = 0$ in \eqref{eq:effStochasticTermsLin},
\eqref{eq:effStochasticTermsResid}, and taking into account that $U^N = \bar{U}^N+R^N$, we conclude
\begin{equation*}%\label{eq:ConsistencyStochTerms}
     \lim_{N \rightarrow \infty}
     \frac{\int_0^T \langle A^{1 + 2 \alpha -\gamma}
       U^N, \sum_{k=1}^{N}  \Phi_{k} dW_{k} \rangle}
     { \int_0^T  | A^{1 + \alpha } U^N |^2 dt } = 0
      \quad \textrm{ in Probability. }
\end{equation*}
Thus, by representation \eqref{eq:RepresentationNuTilde} we have that
$\widetilde{\nu}_{N}$ is weakly consistent estimator of the true parameter $\nu$.
\end{Cor}

We turn next to the `nonlinear terms' appearing in \eqref{eq:RepresentationNuCheck}.
\begin{Prop}\label{thm:ConsistencyInNonlinearTermPlusRates}
 Assume the conditions and notations imposed for Lemma~\ref{thm:ExactEst1} above.
  We suppose that $\delta \in [0, \min\{5/4 - \gamma, \alpha - \gamma + 1\})$ in the
  case (\ref{eq:BCDirchlet}) or that $\delta \in [0, \min\{1/2, \alpha
- \gamma + 1\})$ when we assume (\ref{eq:PeriodicBC}).  Then
\begin{equation}\label{eq:NonlinearConsistencyRate}
   \lim_{N \rightarrow \infty} N^{\delta}
      \frac{\int_0^T \langle A^{1 + 2\alpha} U^N, P_NB(U) \rangle dt}
      {\int_0^T |A^{1+ \alpha } U^N|^2dt} = 0 \quad a.s.
 \end{equation}
\end{Prop}
\begin{proof}
  By the Cauchy-Schwartz inequality,
  \begin{equation}\label{eq:SufficentCond}
    \begin{split}
      \left| \frac{\int_0^T \langle A^{1+ 2\alpha} U^N, P_N B(U) \rangle dt}{\int_0^T |A^{1 +\alpha} U^N|^2dt} \right|
        \leq \left(\frac{\int_0^T |A^{\alpha} P_N B(U)|^2dt}{\int_0^T |A^{1 +\alpha} U^N|^2dt}\right)^{1/2} \, .
    \end{split}
  \end{equation}
  Due to Lemma~\ref{thm:ExactEst1}, \eqref{eq:ExactEst1} and
  \eqref{eq:BlowUpRateLinearSystem} it is therefore sufficient to show that
  \begin{equation}\label{eq:FinalSufficentCondition}
    \lim_{N \rightarrow \infty} \lambda_N^{2(\delta - (\alpha - \gamma + 1))} \int_0^T |A^{\alpha} P_N B(U)|^2dt
        =0 \quad \textrm{a.s.}
  \end{equation}

  We begin with the boundary conditions (\ref{eq:PeriodicBC}) and consider two possibilities
  corresponding to different values of $\alpha$.  First suppose that
  $\alpha < \gamma - 1/2$, so that   $\delta - (\alpha - \gamma +1) < 0$.
  Pick any
  $\beta \in (\max\{\alpha,1/2\}, \gamma -1/2)$.
  Making use of (\ref{eq:HigherOrderEstNonlinearTerm}) and
  then applying
  Theorem~\ref{thm:WellPosedNessandHigherRegularity}, (ii),
  we observe that
  \begin{equation*}%\label{eq:BRegSmallAlpha}
    \begin{split}
      \int_0^T |A^\alpha P_N B(U)|^2dt
      \leq& \int_0^T |A^\beta B(U)|^2dt\\
      \leq& c\int_0^T |A^\beta U|^2 |A^{\beta+1/2} U|^2 dt< \infty
      \quad a.s.
    \end{split}
  \end{equation*}
 and \eqref{eq:FinalSufficentCondition} follows.

   Now suppose that  $\alpha \geq \gamma - 1/2$.
  In this case we pick an element
  $\alpha' \in (\max\{\delta + \gamma - 1, 1/2\} , \gamma - 1/2)$.
  Note that, by assumption, $\delta < 1/2$ so that this
  interval is non-trivial.
  Clearly $\alpha'  < \alpha$ and we apply
  (\ref{eq:decompEstimates}) and again
  \eqref{eq:HigherOrderEstNonlinearTerm} in order to estimate
 \begin{equation}\label{eq:BRegBigAlpha}
    \begin{split}
      \int_0^T |A^\alpha P_N B(U)|^2dt
      & \leq \lambda_N^{2(\alpha - \alpha')} \int_0^T |A^{\alpha'} P_N B(U)|^2dt\\
      & \leq c \lambda_N^{2(\alpha - \alpha')} \int_0^T |A^{\alpha'} U|^2 |A^{\alpha'+1/2} U|^2dt < \infty
      \quad a.s.\\
    \end{split}
  \end{equation}
  As above we find that the quanity on the right hand side is finite
  due to Theorem~\ref{thm:WellPosedNessandHigherRegularity}, (ii).
  Noting that $\delta - (\alpha - \gamma + 1) + \alpha
  - \alpha' < 0$, we infer that  \eqref{eq:FinalSufficentCondition} holds true.

The case of Dirichlet boundary conditions \eqref{eq:BCDirchlet} is addressed in a similar manner.
When $\alpha < 1/4$ we directly apply \eqref{eq:SmallFracOrderEstimate}
to infer \eqref{eq:FinalSufficentCondition}.
When $\alpha \geq 1/4$ we pick any
$\alpha' \in (\delta + \gamma - 1, 1/4)$.  Noting that the conditions on
$\delta$ ensure that this interval is nontrivial and that $\alpha' < \alpha$,
we apply \eqref{eq:decompEstimates} and \eqref{eq:SmallFracOrderEstimate},
in a similar manner to \eqref{eq:BRegBigAlpha} and
infer \eqref{eq:FinalSufficentCondition}
for this case too. The proof is now complete.
\end{proof}

\begin{Cor} In similar manner one can establish the same results as above for $P_NB(U^N)$.
  In particular, for $\delta=0$ we have
 \begin{align}\label{eq:NonlinearConsistencyWORate}
   \lim_{N \rightarrow \infty}
      \frac{\int_0^T \langle A^{1 + 2\alpha} U^N, P_NB(U) \rangle dt}
      {\int_0^T |A^{1+ \alpha } U^N|^2dt} = 0 \quad \textrm{a.s.} \\
      \lim_{N \rightarrow \infty}
      \frac{\int_0^T \langle A^{1 + 2\alpha} U^N, P_NB(U^N) \rangle dt}
      {\int_0^T |A^{1+ \alpha } U^N|^2dt} = 0 \quad \textrm{a.s.}
 \end{align}
Taking into account the above equalities and the representations \eqref{eq:RepresentationNuCheck} and \eqref{eq:RepresentationNuHat} we have
that $\widehat{\nu}_N$ and $\check{\nu}_N$ are consistent estimators of $\nu$.
\end{Cor}

\subsection{Asymptotic Normality}
\label{sec:asymtotic-normality}
We finally address the asymptotic normality of $\tilde{\nu}_N$ and prove
second part of Theorem \ref{thm:MainResultTechincal}. Using the representation
\eqref{eq:RepresentationNuTilde} for $\widetilde{\nu}_N$, Lemma~\ref{thm:ExactEst1}, and
\eqref{eq:ExactEst1} we see that is suffices to establish that
\begin{equation}\label{eq:EstimateForCLT}
    \lim_{N \rightarrow \infty} N
  \frac{\int_0^T \langle A^{1 + 2\alpha - \gamma} \bar{U}^N
      , \sum_k \Phi_k dW_k\rangle}
    {\E\int_0^T |A^{1+ \alpha } \bar{U}^N|^2dt}
    \overset{d}{=} \eta \, ,
 \end{equation}
where $\eta$ is normal random variable with mean zero and variance,
$\frac{2\nu(\alpha-\gamma+1)^2}{\lambda_1T(\alpha-\gamma+1/2)}$
and that
  \begin{equation}\label{eq:ResTermDiesDist}
    \lim_{N \rightarrow \infty} N
  \frac{\int_0^T \langle A^{1 + 2\alpha - \gamma} R^N
      , \sum_k \Phi_k dW_k\rangle}
    {\E\int_0^T |A^{1+ \alpha } \bar{U}^N|^2dt}
    = 0 \quad \textrm{ in Probability.}
  \end{equation}

We establish \eqref{eq:EstimateForCLT} with the aid of Lemma~\ref{thm:MartingaleCLT}.
Let $\sigma_k := \lambda_k^{1 + 2\alpha - \gamma} \bar{u}_k$, and
$\xi_k := \int_0^T \sigma_k^2 \,dt, \ k\geq 1$.
Notice that, due to (\ref{eq:specL2ubark}),
 \begin{equation*}%\label{eq:meanSquareandForthEstSigma}
    \begin{split}
    \E  [ \xi_k ] \sim
   \lambda_k^{2 + 4\alpha - 2\gamma}\lambda_k^{-(1 +2 \gamma)}
   &= \lambda_k^{1 + 4\alpha - 4\gamma}\, ,\\
   \mathrm{Var} [\xi_k ]\sim
   \lambda_k^{4 + 8\alpha - 4\gamma}\lambda_k^{-(3 +4 \gamma)}
   &= \lambda_k^{1 + 8\alpha - 8\gamma}.
   \end{split}
 \end{equation*}
 Define $b_n := \sum_{k =1}^n \E \xi_k$.  Under the
 given assumptions, $1+ 4\alpha - 4 \gamma < -1$, so that
 by \eqref{eq:partialSum} we have that
 $b_n \sim \lambda_n^{2 + 4\alpha - 4\gamma}$.  We infer that
 $b_n$ is increasing and unbounded.  Moreover,
 \begin{equation*}%\label{eq:SummationCondForLNNTowardNormality}
   \sum_{k = 1}^\infty \frac{\mathrm{Var}[\xi_k]}{b_k^2}
   \leq c \sum_{k =1}^\infty k^{-3},
 \end{equation*}
and therefore by LLN, Lemma~\ref{thm:LLN}, we conclude
  \begin{equation*}%\label{eq:SufficentCondNormality}
    \lim_{N \rightarrow \infty}
    \frac{\sum\limits_{k =1}^{N} \xi_k }{\sum\limits_{k =1}^{N} \E \xi_k} =1 \quad \textrm{a.s.}
  \end{equation*}
Consequently, by Lemma~\ref{thm:MartingaleCLT} with $\sigma_k$ defined above, we have
  \begin{equation}\label{eq:MartCLTConcl}
    \lim_{N \rightarrow \infty} \frac
    {\int_0^T \langle A^{1 + 2\alpha - \gamma} \bar{U}^N , \sum_k \Phi_k dW_k\rangle}
      {\left(\E \int_0^T |A^{1+ 2\alpha - \gamma } \bar{U}^N|^2dt \right)^{1/2}} \overset{d}{=} \mathcal{N}(0,1).
  \end{equation}
 Noting that both $1 + \alpha > \gamma$, $1 + 2 \alpha - \gamma > \gamma$
 we may apply (\ref{eq:BlowUpRateLinearSystem}) and infer
 \begin{equation}\label{eq:RatioForNormalDispRate}
   \frac
   {\left(\E \int_0^T |A^{1+ 2\alpha - \gamma } \bar{U}^N|^2dt \right)^{1/2}}
   {\E\int_0^T |A^{1+ \alpha } \bar{U}^N|^2dt}
   \approx \sqrt{\frac{2\nu}{\lambda_1 T}} \cdot  \frac{\alpha-\gamma+1}{\sqrt{\alpha-\gamma+1/2}} \cdot \frac{1}{N}\, .
\end{equation}
By combining (\ref{eq:MartCLTConcl}) and (\ref{eq:RatioForNormalDispRate}) we obtain (\ref{eq:EstimateForCLT}).

The other condition, \eqref{eq:ResTermDiesDist}, follows directly from
Lemma~\ref{thm:ConsistenccyIntheStochasticTerms},
\eqref{eq:effStochasticTermsResid} with the admissable value of $\delta_2 = 1$.
This completes the proof of Theorem~\ref{thm:MainResultTechincal}, (ii).

\begin{Rmk}
Notice that the accuracy of the estimators
$\widetilde{\nu}_{N}$, as measured by
the variance from asymptotic normality, depends on
$\alpha, \ \gamma$ and $T$.
As one would expect, as $T$ gets larger  and more information
is revealed, the quality of the estimator improves (i.e. the variance decreases).
This suggests that one may show that each of the classes of estimators considered
above are consistent in the large time asymptotics regime, i.e. when we fix $N$ and
send $T \rightarrow \infty$.  As with the spectral method we have developed in this
work establishing the long time asymptotics is more complicated in comparison to the
linear case and will be addressed in future work.

Also, we note that $(\alpha-\gamma+1)/\sqrt{\alpha-\gamma+1/2}$,
as a function of $\alpha$ on the domain $\alpha>\gamma-1/2$ reaches its minimum at
$\alpha=\gamma$. Thus, for fixed $T$ and $\gamma$, the smallest
asymptotic variance for the estimator
$\widetilde{\nu}_N$ corresponds to $\alpha=\gamma$.
Observe that when $\alpha=\gamma$
the estimator $\widetilde{\nu}_N$ reduces to
the formal MLE \eqref{eq:EstimatorNuTilde}, which in
some sense is the optimal estimator in this class of estimators.
\end{Rmk}

\begin{Rmk}
We want to emphasize that we still believe, that asymptotic normality properties similar to
\eqref{eq:asymtoticNormalnuHat} also hold true
for the estimators $\check{\nu}_N$. However, for a rigorous
proof one needs to show, for example,
that  Proposition \ref{thm:ConsistencyInNonlinearTermPlusRates} holds
true for some $\delta\geq1$, and
with $P_NB(U)$ replaced by $P_NB(U)-P_NB(U^N)$. Intuitively
it is clear that the difference $P_NB(U)-P_NB(U^N)$ will make the convergence
to zero in \eqref{eq:NonlinearConsistencyRate} faster,
allowing for a larger $\delta$ comparative to those from the
terms $P_NB(U)$ and $P_NB(U^N)$ considered individually.  We further believe
that the quality of these estimators $\check{\nu}_{N}$ may be optimized in terms
of the free parameter $\alpha$.
Although, up to the present time, we remain unable to establish such quantitative results
about the asymptotic normality of the estimators $\check{\nu}_N$,
we  plan to study these questions at least by means of numerical simulations in
forthcoming work.
\end{Rmk}

\section*{Acknowledgments}
We would like to thank the anonymous referees for their
helpful comments and suggestions which improved
greatly the final manuscript.
IC acknowledges support from the National Science Foundation (NSF) grant
DMS-0908099. The work of NGH was partially supported by the NSF under the grants
DMS-1004638, and NSF-DMS-0906440 and by the Research Fund of Indiana University.

\bibliographystyle{plain}

%
%\noindent Igor Cialenco\\ {\footnotesize
%Department of Applied Mathematics\\
%Illinois Institute of Technology\\
%10 West 32nd Str, Bld E1, Room 208\\
%Chicago, IL 60616-3793 \\
%Web: \url{http://math.iit.edu/\~igor}\\
%Email: \url{igor@math.iit.edu} } \\[.3cm]
%Nathan Glatt-Holtz\\ {\footnotesize
%Department of Mathematics\\
%Indiana University\\
%303 Swain East\\
%Bloomington, IN  47405\\
%Web: \url{http://mypage.iu.edu/\~negh/}\\
% Email: \url{negh@indiana.edu}
%}

\end{document}